\documentclass[11pt]{article}
\usepackage{amsmath, amsfonts}

\usepackage{enumerate}% http://ctan.org/pkg/enumerate
\usepackage{color}
\numberwithin{equation}{section}

\def\B1{B_{1/2}}

\def\Bl1{{\lambda_1}}
\def\Box{\hfill\rule{2.5mm}{2.5mm}}
\def\C{{\cal {C}}}

\def\H{{\cal H}}

\def\R{{\mathbb {R}}}

\def\bl1{{\bar \lambda_1}}

\def\build#1_#2^#3{\mathrel{
\mathop{\kern 0pt#1}\limits_{#2}^{#3}}}

\def\h1{\mathop{\rm H^1_{\rm loc,\rm u}}}

\def\l2{\mathop{\rm L^2_{\rm loc,\rm u}}}

\newtheorem{cor}{Corollary}[section]

\newtheorem{lem}[cor]{Lemma}
\newtheorem{prop}[cor]{Proposition}
\newtheorem {rem}  {Remark}
\newtheorem{propo}{Proposition}
\newtheorem{thm}[propo]{Theorem}
\newtheorem{corrol}[propo]{Corollary}
\newtheorem{proposi}[propo]{Proposition}

\begin{document}

\date{}
\title{\bf The blow-up rate for strongly perturbed semilinear wave equations in the conformal case \footnote{The authors are partially supported by the ERC Advanced Grant no.291214, BLOWDISOL during their visit to LAGA, Universit\'{e} Paris 13 in 2013.  }}

\author {M.A. Hamza and O. Saidi\\
%{\it \small Facult\'e des Sciences de Tunis, Universit\'e Tunis El Manar, Tunis, Tunisie}
}
\maketitle
\begin{abstract}
We consider in this work some class of strongly perturbed for the semilinear wave equation with conformal power nonlinearity.\,\,We obtain an optimal estimate for a radial blow-up solution and we have also obtained two less stronger estimates. These results are achieved in three-steps argument by the construction of a Lyapunov functional in similarity variables and the Pohozaev identity derived by multiplying equation $\eqref{p'}$ by $y\partial_{y} w$.
\end{abstract}

\medskip

\section{Introduction}
This paper is devoted to the study of blow-up solutions for the following semilinear
wave equation:
\begin{equation}\label{oy}
\left\{
  \begin{array}{ll}
   \partial_{t}^2 U=\Delta U+|U|^{p_{c}-1}U+f(U)+g(x,t,\nabla U ,\partial_{t} U )\\
  U(0)=U_{0},\partial_{t}U(0)=U_{1},\,\,\,\,\,\,\,\,\,\,\,\,\,\,\,\,\,\,\,\,\,\,\,\\ 
  \end{array}
\right. 
\end{equation}
with conformal power nonlinearity
\begin{equation} \label{65}
p_{c}\equiv 1+\frac{4}{N-1},\,\,{\rm where}\,\,N\geq 2\\
\end{equation}
and $U(t):\,\, x\in {\R^N} \rightarrow U(x,t) \in {\R}$, $U_{0}\in H^{1}_{loc,u}$ and $U_{1} \in L^2_{loc,u}$.
The space $L^2_{loc,u}$ is the set of all $v\in L^2_{loc}$ such that
$$\|v\|_{L^2_{loc,u}}\equiv \sup_{d\in{\R^N}}\Big(\int_{|x-d|<1} |v(x)|^2 dx \Big)^{\frac{1}{2}}<+\infty ,$$
and the space $ H^{1}_{loc,u}=\{v \mid v,|\nabla v|\in L^{2}_{loc,u}\}.$
\\We assume that the functions $f$ and $g$ are $\C^1$, with $f :\R \rightarrow \R$ and $g: \R^{2N+2} \rightarrow \R$ globally lipschitz, satisfying the following conditions:
$$(H_{f})\,\,\,\,\,\,\,\,\,\,\,\,\,\,\,\,|f(v)|\leq M\Big(1+\frac{|v|^{p_{c}}}{\log^{a} (2+v^2)}\Big),\,\,\,\,\,\,{\rm for }\,\,{\rm all} \,v\in \R\,\,{\rm with} \,(M > 0,\,\,\,a> 1),$$
$$(H_{g})\,\,\,\,\,\,\,\,\,|g(x,t,v,z)|\leq  M(1+|v|+|z|), \,\,\,\,\,\,\,\,\,\,\,\,\,{\rm for }\,\,{\rm all} \,x,\,\,v\,\,\in \R^N\,t,z\in \R\,{\rm with }\,(M > 0).$$

\bigskip

The Cauchy problem of equation $\eqref{oy}$ is wellposed in $H^{1}_{loc,u}\times L^{2}_{loc,u}$. This follows from the finite speed of propagation and the wellposdness in $H^{1}\times L^{2}$, valid whenever
$1< p_{c} <1+\frac{4}{N-2}$. It is also known that the existence of blow-up solutions $U(t)$ of $\eqref{oy}$ follows from ODE technics or the energy-based blow-up criterion of Levine \cite{ha} (see
 for example Levine and Todorova \cite{HAG} and Todorova \cite{GT}). More blow-up results can be founded in Caffarelli and Friedman \cite{LA}, \cite{LA1}, Kichenassamy and Littman \cite{SW}, \cite{SW1}, Killip and Visan \cite{RM}. If $U(t)$ is a blow-up solution of
$\eqref{oy}$, we define (see for example Alinhac \cite{sa} and \cite{sa1}) $\Gamma$ as the graph of a function $x\mapsto T(x)$ such that the domain of definition of $U$ (also called the maximal influence domain)
$$D_{U}=\{ (x,t)| t< T(x)\}.$$
The surface $\Gamma$ is called the blow-up graph of $U$.\,\,A point $x_{0} \in \R^N$ is a non-characteristic point if there are:
\begin{equation} \label{17}
 \exists \delta_{0}=\delta_{0}(x_{0})\in (0,1)\,\,{\rm such }\,\,{\rm that} \,\,U\,\, {\rm is}\,\, {\rm defined}\,\, {\rm on}\,\, C_{x_{0},T(x_{0}),\delta_{0}},\\
 \end{equation}
 where
 $$C_{x,t,\delta}=\{(\xi,\tau)\neq (x,t)| 0\leq \tau \leq t- \delta | \xi-x|\}.$$

\bigskip

In the pure power case, equation $\eqref{oy}$ reduces to the semilinear wave equation:
\begin{equation} \label{1.6}
  \partial_{t}^2 U =\Delta U+|U|^{p_{c}-1}U,\,\,\,(x,t)\in \R^N \times [0,T).\\
  \end{equation}
It is interesting to recall that previously Merle and Zaag in \cite{fh3} and \cite{kl} have proved, that if $U$ a solution of $\eqref{1.6}$ with blow-up graph $\Gamma : \{ x\mapsto T(x)\}$ and $x_{0}$ is a non-characteristic point (in the sense $\eqref{17}$) and $1<p\leq p_{c}$, then for all $t\in [\frac{3T(x_{0})}{4},T(x_{0})]$,
 \begin{equation}\label{u}
  0 < \varepsilon_{0}(N,p)\leq (T(x_{0})-t)^{\frac{2}{ p-1}}\frac{\|U(t)\|_{L^{2}(B(x_{0},T(x_{0})-t))}}{(T(x_{0})-t)^{\frac{N}{2}}}
  \end{equation}
  $$+(T(x_{0})-t)^{\frac{2}{ p-1}+1}\Big(\frac{\|\partial_{t} U(t)\|_{L^{2}(B(x_{0},T(x_{0})-t))}}{(T(x_{0})-t)^{\frac{N}{2}}}+\frac{\|\nabla U(t)\|_{L^{2}(B(x_{0},T(x_{0})-t))}}{(T(x_{0})-t)^{\frac{N}{2}}}\Big)\leq K,$$
where the constant $K$ depends only on $N$ and on an upper bound on $T(x_{0})$, $\frac{1}{T(x_{0})}$, $\delta_{0}(x_{0})$ and the initial data in $H^{1}_{loc,u}(\R^N)\times L^{2}_{loc,u}(\R^N)$. 

\bigskip

The unperturbed case \eqref{1.6} is considered in the mathematical community as a lab model for the development of efficient tools for the study of blow-up. Unfortunatly, in more physical situations, the models are often more rich, hence more complicated, with dissipative terms (involving $ \partial_{t} U  $, $\Delta (\partial_{t}U)) $ and other lower order source terms (for example if  $\frac{f(U)}{|U|^{p}}\longrightarrow 0 $ as $U\longrightarrow \infty $ or the case where this term is in the form 
$f=f(x,U)=V(x)|U|^{p}$ where  $V(x)\to 0$ as $x\to 0$). 
 Therefore, it is completely  meaningful for the mathematician to try to extend his methods and results to perturbations of the lab models, since the perturbed models are more encountered in the real-worlds models (see Whitham \cite{W}).

\bigskip

We should recall that in \cite{MH} and \cite{MH1} Hamza and Zaag consider a similar class of perturbed equations, with $(H_{f})$ replaced by a more restrictive conditions: $|f(U)|\leq M(1+|U|^q)$ $q<p,\,\,M > 0$ and they proved a similar result as $\eqref{u}$ valid when the exponent $p$ is subconformal or conformal (i.e. $1<p\leq p_{c}$). Also, when $(H_{f})$ holds Hamza and Saidi \cite{X} have proved a similar result as $\eqref{u}$ valid when $1<p< p_{c}$. However, their methods breaks down in the conformal case (i.e. when $p\equiv p_{c}$), so this case is the subject of this paper.\,\,More precisely, we would rather investigate the growth estimate for $U$ near the space time blow-up graph and extend the result of Hamza and Zaag \cite{MH1} to a stronger class of perturbation, as in \cite{X} when $(H_{f})$ and $(H_{g})$ holds, Hamza and Saidi extend the result of Hamza and Zaag \cite{MH} to a stronger class of perturbation in the subconformal case. We show here that the blow-up rate remains unchanged with perturbations satisfying $(H_{f})$ and $(H_{g})$ and we obtain three different estimates for the blow-up solution, when $a>2$ the estimate of a radial blow-up solution is optimal. However when $1<a<2$, the estimate is not optimal in the radial and general case.  

\bigskip

Among other technics we use in this paper the Pohozaev identity, for that reason we would like to cite some works in which this identity is used. Many results have been established by using the Pohozaev identity, begining with the celebrated result of Pohozaev \cite{sip} where he states that any solution of $\Delta u+f(u)=0$ satisfies an identity, which is known as the Pohozaev identity. The most immediate consequence being the nonexistence of nontrivial bounded solutions for supercritical nonlinearities $f$. The same type of identity is extended by Ros-Oton and Serra \cite{xj} to the semilinear Dirichlet problem namely $(-\Delta )^s u=f(u)$ with $s\in (0,1)$. As an application to the Pohozaev identity the same authors deduce the nonexistence of nontrivial solutions in star-shaped domain for supercritical nonlinearities. We can also mention that among other technics Giga and Kohn in \cite{gk} use the Pohozaev identity to prove that an upper bound on the blow-up solution of the semilinear heat equation is availabale for $1<p<\frac{3N+8}{3N-4}$ or for non-negative initial data with subcritical $p$, then in \cite{gk1} they classify all stationary solutions in self-similar variables. Recently, the same type of identity have been used in the analysis of elliptic PDEs (see \cite{ds} and \cite{sv}).

\bigskip

Let us mention that, the blow-up question for a logarithmic perturbations of pure power nonlinearities is also asked by Nguyen in \cite{VTN} for the following semilinear heat equation:
$$
\left\{
  \begin{array}{ll}
   \partial_{t} U=\Delta U+|U|^{p-1}U+h(U)\\
  U(0)=U_{0}\in L^{\infty},\,\,\,\,\,\,\,\,\,\,\,\,\,\,\,\,\,\,\,\,\,\,\,\\ 
  \end{array}
\right. 
$$
where $p$ is sub-critical nonlinearity (i.e. $1<p\,\,\,{\rm and}\,\,\, p< 1+\frac{4}{N-2}$) and
the function $h$ is in $\C^2(\R\setminus \lbrace 0 \rbrace,\R)$ satisfying
$$(H_{h})\,\,\sum _{j=0}^{2}|z|^j|h^{(j)}(z)|\leq M\Big(1+\frac{|z|^{p}}{\log^{a} (2+z^2)}\Big),\,{\rm or }\,\,h(z)=\xi \frac{|z|^{p-1}z}{\log^{a} (2+z^2)},$$
where $(M > 0,\,\,\,a> 1\,\,{\rm and}\,\,\xi \in \R)$.
\\Nguyen extend the result of Giga, Matsui and Sasayama \cite{GMS} for a logarithmic perturbations of type $(H_{h})$.\,\,He use a Lyapunov functional among other results to derive the blow-up rate.

\bigskip

We would like to mention the remarkable result of Donninger and Sch\"{o}rkhuber in \cite{DS} who proved, in the subconformal range, the stability of the ODE solution $u(t) = \kappa_0 (p)(T - t)^{-\frac{2}{p-1}} $ among all radial solutions, with respect to small perturbations in initial data in the energy topology.\,\,Their approach is based in particular on a good understanding of the spectral properties of the linearized operator in self-similar variables, operator which is not self-adjoint.\,\,Similar results have also been obtained by the same authors \cite{DS1} in the superconformal case and even in the Sobolev supercritical case, (i.e. for any $p>p_{c})$. They extend  to this range the stability result obtained in the subconformal range in \cite{DS}, though they need a topology stronger than needed by the energy. 
Also, in \cite{RBM} Killip, Stoval and Visan are interested to the question of the blow-up rate in the superconformal case  and Sobolev subcritical range (i.e. $N\geq 3$ and $p_{c}<p< p_{s}\equiv \frac{N+2}{N-2}$), there  they  proved that an upper bound is available  for the equation \eqref{1.6} and even the nonlinear Klein-Gordan equation $ \partial_{t}^2 U =\Delta U+|U|^{p-1}U-U$. They  construct a Lyapunov functional in the original variables  by exploiting some  dilation identity unlike in \cite{kl}, \cite{MH}, \cite{MH1} and \cite{X} where  the Lyapunov functional is constructed in  similarity variables. The result of \cite{RBM}  was further refined by Hamza and Zaag in \cite{MH3}.

\bigskip

 Willing to be as exhaustive as possible in our bibliography about the blow-up question for equation $\eqref{1.6}$, we would like to mention some blow-up results in the Sobolev critical range (i.e. $N\geq 3$ and $p= p_{s}\equiv \frac{N+2}{N-2}$), the pure power nonlinearity case $\eqref{1.6}$ has attracted a lot of interest. Many authors addressed the question of obtaining sufficient conditions for scattering and blow-up, through energy estimates, in relation with the ground state (see Kenig and Merle \cite{cef}, Duyckaerts and Merle \cite{TF}). Furthermore, dynamics around the soliton were studied: see Krieger and Shlag \cite{ks}, Krieger, Nakanishi and Shlag \cite{kns} and \cite{kns1}. There are also some remarkable classification theorems by Duyckaerts, Kenig and Merle \cite{TKM}, \cite{TKM1}, \cite{TKM2} and \cite{TKM3}. Analogous results for the critical case of the 
 nonlinear Klein-Gordon equation 
have been proved by Ibrahim, Masmoudi and Nakanishi \cite{IMN} and \cite{IMN1}.

\bigskip

Concerning the blow-up behavior we would want to mention that Donninger, Huang, Krieger and Schlag prove in \cite{Dhk} the existence of so-called "exotic" blow-up solutions when $N=3$, whose blow-up rate oscillates between several pure-power laws.

\medskip

Our method relies on the estimates in similarity variables introduced in \cite{cf} and used in \cite{fh3}, \cite{fh4}, \cite{kl}, \cite{MH}, \cite{MH1} and in \cite{X} (we can see also the radial case treated in \cite{MH2} and \cite{fh8}). More precisely, given $(x_{0},T_{0})$ such that $0< T_{0}\leq T(x_{0})$, we introduce the self-similar change of variables:
\begin{equation}\label{150}
y=\frac{x-x_{0}}{T_{0}-t},\,\,\,\,\,\,s=-\log(T_{0}-t),\,\,w_{x_{0},T_{0}}(y,s)=(T_{0}-t)^\frac{2}{p_{c}-1}U(x,t).\\
\end{equation}
From $\eqref{oy}$, the function $w_{x_{0},T_{0}}$ (we write $w$ for simplicity) satisfies the following
equation for all $y\in B\equiv B(0,1)$ and $s\geq -\log( T_{0})$:
\begin{eqnarray}\label{60}
\partial^2_{s}w &=&div(\nabla w- (y.\nabla w)y)-\frac{2(p_{c}+1)}{(p_{c}-1)^2}w+|w|^{p_{c}-1}w-\frac{p_{c}+3}{p_{c}-1}\partial_{s}w\nonumber \\
&&-2y.\nabla \partial_{s}w+e^{\frac{-2p_{c}s}{p_{c}-1}}f(e^{\frac{2s}{p_{c}-1}}w).
\end{eqnarray}
In the new set of variables $(y,s)$ the behavior of $U$ as $t\rightarrow T_{0}$ is equivalent to the behavior of $w$ as $s\rightarrow +\infty$.

\bigskip

The treatment of the conformal case requires a new idea valid just in the radial case, because the method used in the subconformal case by Hamza and Zaag \cite{MH} and Hamza and Saidi \cite{X} breaks down when $p\equiv p_{c}$, since in the
energy estimates in similarity variables, the perturbations terms are integrated
on the whole unit ball, hence, difficult to control with the dissipation of the non perturbed equation $\eqref{1.6}$, which degenerates to the boundary of the unit ball.
We would like to point out that in the conformal case Hamza and Zaag \cite{MH1} overcame this difficulty via some exponential bound of the $H^{1}\times L^{2}(B)$ norm of the solution. We get here the exponential bound but this estimates are insufficient to conclude our result (see Remark 1 below for explanation).\,\,That obstruction fully justifies our new paper, where we invent a new idea to get our optimal result for a radial blow-up solution of $\eqref{oy}$ when $a>2$.

\bigskip

In what follows we shall fix $f(U)\equiv\frac{|U|^{p_{c}}}{\log^{a} (2+U^2)}$ and $g \equiv 0$, in the equation $\eqref{oy}$. The adaptation to the case $g \not\equiv 0$ is straightforward  from the technics of \cite{MH} and \cite{MH1}. 
\\Now, we announce the following rough estimate: 
\begin{thm}{\bf (Blow-up bounds in the general case).}\label{T1}
Let $a>1$, consider $U$ a solution of $\eqref{oy}$ with blow-up graph $\Gamma : \{ x\mapsto T(x)\}$ and $x_{0}$ is a non-characteristic point (in the sense $\eqref{17}$), then for all $\eta \in (0,1)$, there exists $t_{0}(x_{0})\in [ 0,T(x_{0}))$ such that, for all $T_{0}\in (t_{0}(x_{0}),T(x_{0}))$, for all $s\geq -\log (T_{0}-t_{0}(x_{0}))$ and $y\in B\equiv B(0,1)$, we have
\begin{enumerate}[{\rm(i)}]
\item 
\begin{equation}\label{w2}
\int_{s}^{s+1}\int_{B}(\partial_{s}w(y,\tau))^2\frac{\rho_{\eta}}{1-|y|^2}{\mathrm{d}}y{\mathrm{d}}\tau 
\leq K_{1}e^{\eta \frac{p_{c}+3}{2} s},\\
\end{equation}
\item 
\begin{equation}\label{w3}
\int_{s}^{s+1}\int_{B}|\nabla w(y,\tau)|^2{\mathrm{d}}y{\mathrm{d}}\tau 
+\int_{s}^{s+1}\int_{B}|w(y,\tau)|^{p_{c}+1}{\mathrm{d}}y{\mathrm{d}}\tau\leq K_{1} e^{\eta \frac{p_{c}+3}{2} s},\\
\end{equation}
\end{enumerate}
such that $w=w_{x_{0},T_{0}}$ is defined in $\eqref{150}$ and $\rho_{\eta}$ in $\eqref{oa}$, with 
\\$K_{1}=K_{1}\Big(\eta,T_{0}-t_{0}(x_{0}),\|(U(t_{0}(x_{0})),\partial_{t} U(t_{0}(x_{0}))\|_{H^{1}\times L^2(B(x_{0},\frac{T_{0}-t_{0}(x_{0})}{\delta_{0}(x_{0})}))}\Big)$ and 
\\$\delta_{0}(x_{0})\in (0,1)$ defined in $\eqref{17}$.
\end{thm}

In the original variables Theorem \ref{T1} implies the following:
\begin{corrol} Let $a>1$, consider $U$ a solution of $\eqref{oy}$ with blow-up graph $\Gamma : \{ x\mapsto T(x)\}$ and $x_{0}$ is a non-characteristic point (in the sense $\eqref{17}$), then for all $\eta \in (0,1)$, there exists $t_{0}(x_{0})\in [ 0,T(x_{0}))$ such that for all $t \in [ t_{0}(x_{0}),T(x_{0}))$ we have
\begin{enumerate}[{\rm(i)}]
\item 
$$\int_{T(x_{0})-t}^{T(x_{0})-\frac{t}{2}}\int_{B(x_{0},T(x_{0})-\tau)}|\partial_{t} U(x,\tau)|^{2}{\mathrm{d}}x{\mathrm{d}}\tau  \leq  K_{2}(T(x_{0})-t)^{-\eta \frac{p_{c}+3}{2}},$$
\item 
 $$\int_{T(x_{0})-t}^{T(x_{0})-\frac{t}{2}} \int_{B(x_{0},T(x_{0})-\tau)}\Big(|\nabla U(x,\tau)|^{2}+|U(x,\tau)|^{p_{c}+1}\Big){\mathrm{d}}x{\mathrm{d}}\tau  \leq  K_{2}(T(x_{0})-t)^{-\eta \frac{p_{c}+3}{2}}.$$
\end{enumerate}
\end{corrol}

\bigskip

Note that, Hamza and Zaag \cite{MH1} exploit equation $\eqref{01}$ with some functional associated to obtain an exponentially estimate to the blow-up solution, which induce with a natural Lyapunov functional for equation $\eqref{60}$ to deduce their optimal result.
In this work, it is not the case because our perturbation is stronger than the one in Hamza and Zaag \cite{MH1}.\,\,To overcame this difficulty, we assume that $U$ is a radial blow-up solution of equation $\eqref{oy}$ and we insert between the exponential estimate and the optimal estimate for a radial blow-up solution when $a>2$ a polynomial estimate for a radial blow-up solution when $a>1$ obtained by transforming equation $\eqref{60}$ to equation $\eqref{0301}$ and by exploiting some functional associated with the weight $\rho_{s^{-b}}$ which will be defined later in $\eqref{phi}$, just we need to insist that the weight depends on time.

\begin{rem}
{\rm The result of Theorem 1 is similar to the one obtained in Hamza and Zaag \cite{MH1}, (see Proposition 2.1 page 201). Unfortunately, when ($H_{f}$) holds we can not conclude our optimal result for a radial blow-up solution when $a>2$ because our perturbation are polynomially smaller in time, so the exponential bound is not sufficient to conclude. The idea is to take a weight $\rho_{s^{-b}}$ which will be defined later in $\eqref{phi}$ and we rewrite equation $\eqref{60}$ in the radial case (see $\eqref{p'}$) and also in another form (see $\eqref{0301}$) to obtain a polynomially bound to the blow-up solution.\,\,We need just to explain the difficulty in the general case, the dependance of $\rho_{s^{-b}}$ on time give birth to a new terms of type $-\frac{b}{2s^{b+1}}\int_{B}\Big(\frac{1}{2}(\partial_{s}w)^2+\frac{p_{c}+1}{(p_{c}-1)^2}w^2-\frac{|w|^{p_{c}+1}}{p_{c}+1}+e^{-\frac{2(p_{c}+1)s}{p_{c}-1}}F(e^{\frac{2s}{p_{c}-1}}w)\Big)\log(1-|y|^2)\rho_{s^{-b}}{\mathrm{d}}y$ and $-\frac{b}{2s^{b+1}}\int_{B}\Big(|\nabla w|^2-(y.\nabla w)^2\Big)\log(1-|y|^2)\rho_{s^{-b}}{\mathrm{d}}y$, by some technics we control all these terms except the last one is problematic and following this change of variables $y=r\omega $ with $r=|y|$ and $\omega =\frac{y}{|y|}$ we can see that it is more difficult to handle the tangential part of $\nabla w$, of course under radial symmetry where this term vanishes (see Section 3 for details).\,\,For that reason, we restrict ourselves from now on to radially symmetric data, where we obtain our optimal result (see Theorem 6 below).}
\end{rem}

\medskip

Before entering into the details of our second result and in what follows, we take $U_{0}$, $U_{1}$ a radial initial data and the function $g=g(|x|,t,\nabla U.\frac{x}{|x|},\partial_{t} U )$ in $\eqref{oy}$.
We are going to announce our second result for $U$ a radial blow-up solution of $\eqref{oy}$.
\\Since $U$ is radial, we introduce
$$u(r,t)=U(x,t)\,\,{\rm if}\,\,r=|x|,$$
and rewrite $\eqref{oy}$ as
\begin{equation}\label{y}
\left\{
  \begin{array}{ll}
   \partial_{t}^2 u=\partial_{r}^2 u+\frac{N-1}{r}\partial_{r} u+|u|^{p_{c}-1}u+f(u)+g(r,t,\partial_{r}u,\partial_{t} u ),\\
  \partial_{r}u(0,t)=0,\\
  u(r,0)=u_{0}(r)\,\,{\rm and}\,\,\partial_{t}u(r,0)=u_{1}(r),\,\,\,\,\,\,\,\,\,\,\,\,\,\,\,\,\,\,\,\,\,\,\,\\ 
  \end{array}
\right. 
\end{equation}
where $u(t):\,\, r\in \R^+ \mapsto u(r,t) \in {\R}$. 
We use the arguments of Hamza and Zaag \cite{MH2} and also used by Merle and Zaag \cite{fh8}, for a perturbed semilinear wave equation in the radial case, what we brought a new idea when we use this argument is that we include the origin to the set of the non-characteristic point, which is the most important case because in this previous work they are unable to get this result in the origin. Note that, if we are far from the origin we can exploit the technics used by Merle and Zaag \cite{fh8} and Hamza and Zaag \cite{MH2}, where they exclude the origin which bring a singular term $\frac{N-1}{r}\partial_{r} u$ to $\eqref{y}$ combined with the technics used by Hamza and Saidi \cite{X}. Let us briefly explain how we treat the conformal case for a radial blow-up solution when we are outside the origin (i.e. $r_{0}>0$) and $a>1$.
We can see following the change of variable (\ref{150}) that $w$ is a solution of 
\begin{eqnarray}\label{fin}
\partial^2_{s}w &=&\frac{1}{\rho} \partial_{y}(\rho(1-y^2)\partial_{y}w)-\frac{2(p_{c}+1)}{(p_{c}-1)^2}w+|w|^{p_{c}-1}w-\frac{p_{c}+3}{p_{c}-1}\partial_{s}w\nonumber \\
&&+e^{-s}\frac{N-1}{r_{0}+ye^{-s}}\partial_{y}w-2y \partial^2_{y,s}w+e^{\frac{-2p_{c}s}{p_{c}-1}}f(e^{\frac{2s}{p_{c}-1}}w),
\end{eqnarray}
where 
$$\rho(y)=(1-y^2)^{\frac{2}{p_{c}-1}}.$$
We can remark that if $s$ is large and with the fact that $|y|<1$, we have:
$$e^{-s}\Big|\frac{N-1}{r_{0}+ye^{-s}}\partial_{y}w\Big|\leq \frac{2(N-1)}{r_{0}}e^{-s}|\partial_{y}w|,$$
we can see the last bound as a perturbation when $s$ is large enough.
According to Hamza and Saidi \cite{X} we introduce for $\eqref{fin}$ the following Lyapunov functional:
\begin{equation}
H_{r_{0}>0}(w(s),s)=\exp\Big(\frac{p_{c}+3}{(a-1)s^{\frac{a-1}{2}}}\Big) E_{r_{0}>0}(w(s),s)+\theta e^{-\frac{p_{c}+1}{p_{c}-1}s},\\
\end{equation}
where $\theta$ is a large constant.
\begin{eqnarray*}
E_{r_{0}>0}(w(s),s)\hspace{-0,3cm}&=&\hspace{-0,4cm}\int_{-1}^{1}\Big(\frac{1}{2}(\partial_{s}w)^2+\frac{1}{2}(\partial_{y}w)^2(1-y^2)+\frac{p_{c}+1}{(p_{c}-1)^2}w^2-\frac{1}{p_{c}+1}|w|^{p_{c}+1}\Big)\rho  {\mathrm{d}}y\\
&-&e^{\frac{-2(p_{c}+1)s}{p_{c}-1}} \int_{-1}^{1}F(e^{\frac{2s}{p_{c}-1}}w)\rho dy-\frac{1}{s^{b}}\int_{-1}^{1} w\partial_{s}w \rho dy.
\end{eqnarray*}
In other word we can solve the problem of the conformal case for a radial blow-up solution when we are outside origin in one time, it means that we don't need to go through the exponential and the polynomial bound of Theorem 1 and Theorem 4 treated in this paper.
\\Our interest now to derive the blow-up rate for a radial blow-up solution of $\eqref{oy}$ in the origin. Let $r_{0}=0$ a non-characteristic point (here start the novelties) and if we write for any $T_{0}$ such that $0<T_{0}\leq T(0)$ equation $\eqref{60}$ in the radial case, we obtain the following:
\begin{equation}\label{0420}
w_{0}(y,s)=(T_{0}-t)^\frac{2}{p_{c}-1}u(r,t),\,\,\,\,\,\,y=\frac{r}{T_{0}-t},\,\,\,\,\,\,s=-\log(T_{0}-t).\\
\end{equation}
The function $w=w_{0}$ satisfies the following equation for all $y\in (0,1)$ and 
\\$s\geq -\log(T_{0})$
\begin{eqnarray}\label{p'}
\partial^2_{s}w &=&\frac{1}{y^{N-1}\rho_{\eta}} \partial_{y}(y^{N-1}\rho_{\eta}(1-y^2)\partial_{y}w)+2\eta y\partial_{y}w-\frac{2(p_{c}+1)}{(p_{c}-1)^2}w\nonumber \\
&&+|w|^{p_{c}-1}w-\frac{p_{c}+3}{p_{c}-1}\partial_{s}w-2y \partial^2_{y,s}w+e^{\frac{-2p_{c}s}{p_{c}-1}}f(e^{\frac{2s}{p_{c}-1}}w),
\end{eqnarray}
where $\rho_{\eta}=(1-y^2)^{\eta}.$
We construct a new functional $N_{\eta}(w(s))$ defined in $\eqref{L}$ obtained by multiplying equation $\eqref{p'}$ by $y\partial_{y} w$ to write the Pohozaev identity which is crucial to deduce Proposition 3, this idea is effective just in the radial case. Thanks to this new functional we obtain the following proposition: 
\begin{proposi}\label{p}
For all $\eta \in (0,1)$, there exists $S_{6}\geq 1$ such that for all $s\geq \max(S_{6},s_{0})$ we have the following inequality
\begin{equation}\label{w1}
\int_{s}^{s+1}\int_{0}^{1}|w|^{p_{c}+1}\frac{\rho_{\eta}}{1-y^2}y^{N-1} {\mathrm{d}}y{\mathrm{d}}\tau \leq  Ce^{\frac{\eta(p_{c}+3)s}{2}}.\\
\end{equation}
\end{proposi}

\bigskip

We use now Theorem 1 in a clever way and exceptionally in the neighborhood of the edge of the unit ball combined with Proposition 3 to obtain our second result (see Theorem 4 below). 
\\Let us state our second result which is crucial to deduce the main goal of this paper (see Theorem 6 below).
\begin{thm}{\bf (A polynomially estimate).}
Let $a>1$, consider $u$ a solution of equation $\eqref{y}$ with blow-up graph $\Gamma : \{ r\mapsto T(r)\}$, then there exists $t_{1}(0)\in [ 0,T(0))$ such that, for all $T_{0}\in (t_{1}(0),T(0))$, for all $s \geq -\log(T_{0}-t_{1}(0))$ and $y\in (0,1)$ we have for all $b\in (1,a)$
\begin{enumerate}[{\rm(i)}]
\item 
\begin{equation}\label{w4}
\int_{s}^{s+1}\int_{0}^{1}(\partial_{s}w(y,\tau))^2\frac{\varphi (y,\tau)}{1-y^2}{\mathrm{d}}y{\mathrm{d}}\tau  \leq K_{3}s^{b},\\
\end{equation}
\item 
\begin{equation}\label{w5}
\int_{s}^{s+1}\int_{0}^{1}(\partial_{y}w(y,\tau))^2(1-y^2)\varphi (y,\tau){\mathrm{d}}y{\mathrm{d}}\tau  \leq K_{3}s^{b},\\
\end{equation}
\item 
\begin{equation}\label{w6}
\int_{s}^{s+1}\int_{0}^{1}|w(y,\tau)|^{p_{c}+1}y^{N-1} {\mathrm{d}}y{\mathrm{d}}\tau \leq K_{3}s^{b},\\
\end{equation}
\end{enumerate}
where $K_{3}=K_{3}\Big(K_{1},T_{0},\|(u(t_{1}(0)),\partial_{t} u(t_{1}(0))\|_{H^1\times L^2((-\frac{T_{0}-t_{1}(0)}{\delta_{0}(0)},\frac{T_{0}-t_{1}(0)}{\delta_{0}(0)}))}\Big)$, $\delta_{0}(0)$ is defined in $\eqref{17}$ and 
\begin{equation}\label{phi}
\varphi (y,s)=y^{N-1}\rho_{s^{-b}}\,\,\,\,\,{\rm with}\,\,\,\,\rho_{s^{-b}}=(1-y^2)^{s^{-b}}.\\
\end{equation}
\end{thm}

\bigskip

\begin{rem}
{\rm Let us remark that in $(i)$ and $(ii)$ of Theorem 4 we can not control the 
time average of the $L^2$ norm of $\partial_{s}w$ and $\partial_{y}w$
blow-up solution until the edge of the unit ball.\,\,We get this estimates by using a classic Lyapunov functional obtained by multiplying equation $\eqref{0301}$ by $\partial_{s}w$. However, in $(iii)$ we can control the time average of the $L^{p_{c}+1}$ norm of $w$ until the edge of the unit ball. 
Let us mention that, the following estimate:
$$\int_{s}^{s+1}\int_{0}^{1}|w(y,\tau)|^{p_{c}+1}\varphi (y,\tau) {\mathrm{d}}y{\mathrm{d}}\tau \leq K_{3}s^{b},$$
ensue from Proposition 3.1 which helps us with Proposition 3 where we use essentially the Pohozaev identity to get the estimate $(iii)$ of Theorem 4 as desired.}
\end{rem}

We write now Theorem 4 as the original variables in the following corollary.
\begin{corrol} Let $a>1$, consider $u$ a solution of equation $\eqref{y}$ with blow-up graph $\Gamma : \{ r\mapsto T(r)\}$, then there exists $t_{1}(0)\in [ 0,T(0))$ such that for all 
\\$t\in [ t_{1}(0),T(0))$, we have for all $b\in (1,a)$
\begin{enumerate}[{\rm(i)}]
\item $$\int_{T(0)-t}^{T(0)-\frac{t}{2}} \int_{T(0)-\tau}^{T(0)+\tau}\Big(|\partial_{t} u(r,\tau)|^{2}+|\partial_{r}u(r,\tau)|^{2}\Big){\mathrm{d}}r{\mathrm{d}}\tau\leq K_{4}(-\log(T(0)-t))^{b},$$
 \item  
 $$\int_{T(0)-t}^{T(0)-\frac{t}{2}} \int_{T(0)-\tau}^{T(0)+\tau}|u(r,\tau)|^{p_{c}+1}{\mathrm{d}}r{\mathrm{d}}\tau \leq K_{4}(-\log(T(0)-t))^{b}.$$
 \end{enumerate}
\end{corrol}

\bigskip

Now via Theorem 4 we are in position to announce our main result in the following theorem:
\begin{thm} {\bf (Optimal blow-up rate in the radial case)}
Let $a>2$, consider $u$ a solution of equation $\eqref{y}$ with blow-up graph $\Gamma : \{ r\mapsto T(r)\}$, then there exists $\varepsilon_{0}> 0$ and $\hat{S_{2}}> 0$ such that, 
for all  $s \geq \hat{s_{2}}(0)=\max(\hat{S_{2}},-\log(\frac{T(0)}{4}))$

\begin{equation}\label{w7}
0 < \varepsilon_{0} \leq \|w_{0,T(0)}(s)\|_{H^{1}((0,1))} + \|\partial_{s}w_{0,T(0)}(s)\|_{L^{2}((0,1))}\leq K_{5},\\
\end{equation}
where $K_{5}=K_{5}\Big(K_{3},\hat{s_{2}}(0),\parallel(u(t_{2}(0)),\partial_{t}u(t_{2}(0))\parallel_{H^{1}\times L^{2}((-\frac{e^{-\hat{s_{2}}(0)}}{\delta_{0}(0)},\frac{e^{-\hat{s_{2}}(0)}}{\delta_{0}(0)}))}\Big)$ and 
\\$\delta_{0}(0)\in (0,1)$ defined in $\eqref{17}$.
\end{thm}

\medskip

As for Theorems 1 and 2 we translate Theorem 6 in the original variables, our goal becomes the following corollary:

\newpage

\begin{corrol} Let $a>2$, consider $u$ a solution of equation $\eqref{y}$ with blow-up graph $\Gamma : \{ r\mapsto T(r)\}$, then there exists $\varepsilon_{0}> 0$, such that for all $t \in [ t_{2}(0),T(0))$ with 
\\$t_{2}(0)=T(0)-e^{-\hat{s_{2}}(0)}$, we have 
 \begin{eqnarray*}
  0 < \varepsilon_{0}(N)&\leq &(T(0)-t)^{\frac{2}{ p_{c}-1}}\frac{\|u(t)\|_{L^{2}((t-T(0),T(0)-t))}}{(T(0)-t)^{\frac{N}{2}}}\nonumber\\
  &&+(T(0)-t)^{\frac{2}{ p_{c}-1}+1}\Big(\frac{\|\partial_{t} u(t)\|_{L^{2}((t-T(0),T(0)-t))}}{(T(0)-t)^{\frac{N}{2}}}\\
  &&+
 \frac{\|\partial_{r} u(t)\|_{L^{2}((t-T(0),T(0)-t))}}{(T(0)-t)^{\frac{N}{2}}}\Big)\leq K_{6}.
  \end{eqnarray*}
\end{corrol}

\medskip

\begin{rem} 
{\rm 
 In a series of papers \cite {fh3}, \cite {kl}, \cite {fh4}, \cite {fh}, \cite {fh1} and \cite {fh6}, Merle and Zaag give a full picture of the blow-up for the solutions of $\eqref{y}$ in one space dimension when $(f,g)\equiv (0,0)$. Moreover, the resulting facts of all this papers are extended by Hamza and Zaag for a perturbed semilinear wave equation in one space dimension or in dimension $N\geq 2$ in \cite{MH} and \cite{MH1} and for a radial blow-up solution outside origin in \cite{MH2}. Finally the result of Hamza and Zaag \cite{MH} is extended by Hamza and Saidi \cite{X} for a strongly perturbed semilinear wave equation in the sub-conformal case. As a matter of fact, our focal interest is in studying the conformal case. 
}
\end{rem}

\paragraph{Layout of the paper.} This paper is organized as follow: Section 2 is devoted to obtain a rough control on space-time of the solution $w$.\,\,Based upon this result, in Section 3 we will prove that the exponential bound obtained in the general case turns into a polynomial bound in the radial case when $a>1$. To do that, we construct a Lyapunov functional in similarity variables and a new functional $N_{\eta}(w(s))$ obtained by multiplying equation $\eqref{p'}$ by $y\partial_{y} w$.\,\,Furthermore, the new functional $N_{\eta}(w(s))$ allows us to control the blow-up solution until the edge of the unit ball. In our case, according to the simple fact that our weight $\rho_{s^{-b}}$ defined in $\eqref{phi}$ depends on time, we can easily notice that, compared to the previous work for example Hamza and Zaag \cite{MH}, \cite{MH1}, Hamza and Saidi \cite{X} and Merle and Zaag \cite{kl} the derivative in time give birth to a novel terms which was already controlled, eventually, we conclude Theorem 4. Finally, in Section 4, according to all this results, we use a technics similar to the one used by Hamza and Zaag \cite{MH1} to conclude our optimal result in the radial case when $a>2$, which is the main goal of this paper.

\bigskip

We mention that $C$ will be used to denote a constant that's depends on $N$, $a$ and $M$ which may vary from line to line.\,\,In the whole paper we assume that $\eqref{65}$ holds and we denote by 
\begin{equation}\label{90}
 F(u)=\int_{0}^{u}f(v){\mathrm{d}}v.
 \end{equation}
\section{Proof of Theorem 1}
The proof follows the same pattern as the perturbed case considered by Hamza and Zaag \cite{MH1}, the unique difference lays in the treatment of the perturbed term. We handle this term and we obtain in this section an exponentially growing bound on time averages of the $H^{1}\times L^2(B)$ norm of $(w,\partial_{s} w)$. Consider $U$ a solution of $\eqref{oy}$ with blow-up graph $\Gamma : \{ x\mapsto T(x)\}$ and $x_{0}$ is a non-characteristic point, the aim of this section is to prove Theorem 1.
\subsection{A Lyapunov functional for equation $\eqref{01}$}
Consider $T_{0}\in (0,T(x_{0})] $, then we write $w$ instead of $w_{x_{0},T_{0}}$ defined in $\eqref{150}$. 
\\Let $\eta \in (0,1)$ and we write equation $\eqref{60}$ satisfied by $w$ in the following form 
\begin{eqnarray}\label{01}
\partial^2_{s}w&=&\frac{1}{\rho_{\eta}} div(\rho_{\eta}\nabla w-\rho_{\eta} (y.\nabla w)y)+2\eta(y.\nabla w)-\frac{2(p_{c}+1)}{(p_{c}-1)^2}w\nonumber \\
&&+|w|^{p_{c}-1}w-\frac{p_{c}+3}{p_{c}-1}\partial_{s}w-2y.\nabla \partial_{s}w+
e^{\frac{-2p_{c}s}{p_{c}-1}}f(e^{\frac{2s}{p_{c}-1}}w),
\end{eqnarray}
where 
\begin{equation}\label{oa}
 \rho_{\eta}=(1-|y|^2)^{\eta}.
\end{equation}
The equation $\eqref{01}$ will be studied in the space $\H$
$$\H =\{(w_{1},w_{2})| \int_{B} \Big(w_{2}^2+|\nabla w_{1}|^2(1-|y|^2)+w_{1}^2\Big)\rho {\mathrm{d}}y< +\infty\}. $$
To control the norm of $(w(s),\partial_{s} w(s) ) \in \H$, we first introduce the following functionals
\begin{eqnarray}\label{10}
\hspace{-0,5cm}E_{\eta}(w(s),s)\hspace{-0,3cm}&=&\hspace{-0,5cm}\int_{ B}\Big(\frac{1}{2}(\partial_{s}w)^2+\frac{1}{2}|\nabla w|^2-\frac{1}{2}(y.\nabla w)^2+\frac{p_{c}+1}{(p_{c}-1)^2}w^2-\frac{1}{p_{c}+1}|w|^{p_{c}+1}\Big)\rho_{\eta} {\mathrm{d}}y\nonumber\\
&&-e^{\frac{-2(p_{c}+1)s}{p_{c}-1}} \int_{ B} F(e^{\frac{2s}{p_{c}-1}}w)\rho_{\eta} {\mathrm{d}}y,\nonumber\\
J_{\eta}(w(s))&=&-\eta\int_{B} w\partial_{s}w \rho_{\eta} dy+\frac{N\eta}{2}\int_{ B}w^2\rho_{\eta} {\mathrm{d}}y ,\\
H_{\eta}(w(s),s)&=&E_{\eta}(w(s),s)+J_{\eta}(w(s)),\nonumber \\
G_{\eta}(w(s),s)&=&H_{\eta}(w(s),s)e^{\frac{-\eta(p_{c}+3)s}{2}}+\theta e^{\frac{-\eta(p_{c}+3)s}{2}},\nonumber
\end{eqnarray}
where $\theta=\theta(\eta)$ is a sufficiently large constant that will be determined later. 
In this subsection we prove that $G_{\eta}(w(s),s)$ is decreasing in time, which will give the rough (i.e exponentially fast) estimate for $E_{\eta}(w(s),s)$ and the time average of the $\|(w,\partial_{s} w)\|_{H^{1}(B)\times L^2(B)}$. More precisely we are going to prove the following proposition

\begin{prop}
For all $\eta \in (0,1)$, there exists $S_{0}\geq 1$ and $\lambda_{0}>0$ such that $G_{\eta}(w(s),s)$ defined in $\eqref{10}$ satisfies for all $s\geq \max (s_{0},S_{0})$,
\begin{eqnarray*} 
G_{\eta}(w(s+1),s+1)-G_{\eta}(w(s),s)&\leq &-\lambda_{0}\int_{s}^{s+1}e^{\frac{-\eta(p_{c}+3)\tau}{2}}\int_{B}(\partial_{s} w)^2\frac{\rho_{\eta}}{1-|y|^2}{\mathrm{d}}y{\mathrm{d}}\tau\nonumber \\
&&-\lambda_{0}\int_{s}^{s+1}e^{\frac{-\eta(p_{c}+3)\tau}{2}}\int_{B}|w|^{p_{c}+1}\rho_{\eta} {\mathrm{d}}y{\mathrm{d}}\tau\nonumber \\
&&-\lambda_{0}\int_{s}^{s+1}e^{\frac{-\eta(p_{c}+3)\tau}{2}}\int_{B}|\nabla w|^2(1-|y|^2)\rho_{\eta} {\mathrm{d}}y{\mathrm{d}}\tau ,
\end{eqnarray*} 
where $w=w_{x_{0},T_{0}}$ is defined in $\eqref{150}$.
Moreover, there exists $S_{1}\geq S_{0}$ such that for all $s\geq max(s_{0},S_{1})$ $G_{\eta}(w(s),s)\geq 0$.
 \end{prop}

\bigskip

Now we state two lemmas which are crucial for the proof. We begin with bounding the time derivative of $E_{\eta}(w(s),s)$ in the following lemma.
\begin{lem}
For all $\eta \in (0,1)$, we have for all $s\geq \max (s_{0},1)$,
 \begin{equation}\label{4}
\frac{d}{ds}(E_{\eta}(w(s),s))= -2\eta\int_{B}(\partial_{s} w)^2\frac{|y|^{2}\rho_{\eta}}{1-|y|^2}{\mathrm{d}}y +2\eta\int_{B}(\partial_{s} w) (y.\nabla w)\rho_{\eta}{\mathrm{d}}y+\Sigma_{0}(s),\\
\end{equation}
where $\Sigma_{0}(s)$ satisfies
    \begin{equation}\label{5}
\Sigma_{0}(s)\leq  \frac{C}{s^{a}}\int_{B}|w|^{p_{c}+1}\rho_{\eta} {\mathrm{d}}y+
Ce^{-\frac{p_{c}+1}{p_{c}-1}s}.
\end{equation}
\end{lem}
\bigskip
{\it Proof}: Multiplying $\eqref{01}$ by $\rho_{\eta}\partial_{s} w$ and integrating over the ball $B$, we obtain $\eqref{4}$ with 
 \begin{equation}\label{ao}
\Sigma_{0}(s)=\frac{2(p_{c}+1)}{p_{c}-1}e^{\frac{-2(p_{c}+1)s}{p_{c}-1}}\int_{B} F(e^{\frac{2s}{p_{c}-1}}w)\rho_{\eta} {\mathrm{d}}y-\frac{2e^{\frac{-2p_{c}s}{p_{c}-1}}}{p_{c}-1}\int_{B}f(e^{\frac{2s}{p_{c}-1}}w)w\rho_{\eta} {\mathrm{d}}y.\\
\end{equation}
Clearly the function $F$ defined in $\eqref{90}$ satisfies the following estimate: 
\begin{equation}\label{104}
|F(x)|+|xf(x)|\leq C\Big(1+\frac{ |x|^{p_{c}+1}}{\log^a(2+x^2)}\Big).\\
\end{equation}
In order to prove $\eqref{5}$, we divide the unit ball $B$ into two parts
 $$A_{1}(s)=\{y \in B\,\,|\,\, w^2(y,s)\leq  e^{\frac{-2s}{p_{c}-1}}\}\,\,{\rm and }\,\,A_{2}(s)=\{y \in B\,\,|\,\, w^2(y,s)>  e^{\frac{-2s}{p_{c}-1}}\}.$$
On the one hand, ${\rm if }\,\,y \in A_{1}(s)$ we have
\begin{equation}\label{93}
\int_{A_{1}(s)}\frac{|w|^{p_{c}+1}}{\log^a(2+e^{\frac{4s}{p_{c}-1}}w^2)}\rho_{\eta} {\mathrm{d}}y \leq \frac{e^{-\frac{p_{c}+1}{p_{c}-1}s}}{\log^a(2)}\int_{A_{1}(s)}\rho_{\eta} {\mathrm{d}}y\leq C e^{-\frac{p_{c}+1}{p_{c}-1}s}. \\
\end{equation}
On the other hand, if $y\in A_{2}(s)$ we have $\log(2+e^{\frac{4s}{p_{c}-1}}w^2)>\frac{2s}{p_{c}-1},$ we obtain for all $s\geq \max (s_{0},1)$
\begin{equation}\label{102}
\int_{A_{2}(s)}\frac{|w|^{p_{c}+1}}{\log^a(2+e^{\frac{4s}{p_{c}-1}}w^2)}\rho_{\eta}{\mathrm{d}}y\leq \frac{C}{s^a}\int_{B}|w|^{p_{c}+1} \rho_{\eta} {\mathrm{d}}y.\\
\end{equation}
To conclude, it suffices to combine $\eqref{93}$ and $\eqref{102}$, then write
\begin{equation}\label{200}
\Sigma_{0}(s) \leq Ce^{-\frac{p_{c}+1}{p_{c}-1}s}+\frac{C}{s^a}\int_{B}|w|^{p_{c}+1} \rho_{\eta} {\mathrm{d}}y,\\
\end{equation}
which ends the proof of Lemma 2.2.
 \Box
       
 \medskip

 We are now going to prove the following estimate for the functional $J_{\eta}(w(s))$.
\begin{lem}
For all $\eta \in (0,1)$, we have for all $s\geq \max (s_{0},1)$,
 \begin{eqnarray}\label{6}
\frac{d}{ds}(J_{\eta}(w(s))&\leq &\frac{32\eta}{p_{c}+15}\int_{B}(\partial_{s} w)^2\frac{\rho_{\eta}}{1-|y|^2}{\mathrm{d}}y +\frac{\eta(p_{c}+3)}{2}H_{\eta}(w(s),s)\\
&&-2\eta\int_{B}\partial_{s} w (y.\nabla w)\rho_{\eta}{\mathrm{d}}y-\frac{\eta(p_{c}-1)}{8}\int_{B}|\nabla w|^2(1-|y|^2)\rho_{\eta} {\mathrm{d}}y \nonumber \\
&&-\frac{\eta(p_{c}+15)}{8}\int_{B}(\partial_{s} w)^2\rho_{\eta}{\mathrm{d}}y-\frac{\eta(p_{c}-1)}{2(p_{c}+1)}\int_{B}|w|^{p_{c}+1} \rho_{\eta}{\mathrm{d}}y+ \Sigma_{1}(s),\nonumber
\end{eqnarray}
where $\Sigma_{1}(s)$ satisfies
    \begin{equation}\label{7}
\Sigma_{1}(s)\leq  \frac{C}{s^{a}}\int_{B}|w|^{p_{c}+1}\rho_{\eta} {\mathrm{d}}y+C\int_{B}w^2\rho_{\eta} {\mathrm{d}}y+Ce^{-\frac{p_{c}+1}{p_{c}-1}s}.
\end{equation}
\end{lem}
{\it Proof}: Note that $J_{\eta}(w(s))$ is a differentiable function, by using equation $\eqref{01}$ and integrating by part, for all $s\geq \max (s_{0},1)$ we have
\begin{eqnarray}\label{20}
\frac{d}{ds}(J_{\eta}(w(s)))\hspace{-0,3cm}&=&\hspace{-0,3cm}-\eta\int_{B}(\partial_{s} w)^2\rho_{\eta}{\mathrm{d}}y+\eta\int_{B}(|\nabla w|^2-(y.\nabla w)^2)\rho_{\eta} {\mathrm{d}}y-2\eta\int_{B} \partial_{s}w(y.\nabla w)\rho_{\eta}{\mathrm{d}}y\nonumber \\
&&-\eta\int_{B}|w|^{p_{c}+1}\rho_{\eta}{\mathrm{d}}y+\underbrace{4\eta^2\int_{B}w \partial_{s}w\frac{|y|^2\rho_{\eta}}{1-|y|^2}{\mathrm{d}}y-2\eta^3\int_{B}w^2\frac{|y|^2\rho_{\eta}}{1-|y|^2}{\mathrm{d}}y}_{\Sigma_{1}^{1}(s)}\nonumber \\
&&+\eta\Big(N\eta+\frac{2p_{c}+2}{(p_{c}-1)^2}\Big)\int_{B}w^2\rho_{\eta} {\mathrm{d}}y-\eta e^{-\frac{2p_{c}s}{p_{c}-1}}\int_{B}wf(e^{\frac{2s}{p_{c}-1}}w)\rho_{\eta}{\mathrm{d}}y. 
\end{eqnarray}
By combining $\eqref{10}$ and $\eqref{20}$ and some straightforward computations, we obtain 
\begin{eqnarray}\label{12}
\frac{d}{ds}(J_{\eta}(w(s)))\hspace{-0,3cm}&= &\hspace{-0,3cm}-2\eta\int_{B} \partial_{s}w(y.\nabla w)\rho_{\eta}{\mathrm{d}}y+\frac{\eta(p_{c}+3)}{2}H_{\eta}(w(s),s)\nonumber \\
&&-\frac{\eta(p_{c}+7)}{4}\int_{B}(\partial_{s} w)^2\rho_{\eta}{\mathrm{d}}y-\eta\Big(\frac{p_{c}+1}{2(p_{c}-1)}+\frac{N\eta(p-1)}{4}\Big)\int_{B}w^2\rho_{\eta} {\mathrm{d}}y\nonumber \\
&&-\frac{\eta(p_{c}-1)}{2(p_{c}+1)}\int_{B}|w|^{p_{c}+1} \rho_{\eta}{\mathrm{d}}y-\frac{\eta(p_{c}-1)}{4}\int_{B}(|\nabla w|^2-(y.\nabla w)^2)\rho_{\eta} {\mathrm{d}}y\nonumber \\
&&-\eta e^{-\frac{2p_{c}s}{p_{c}-1}}\int_{B}wf(e^{\frac{2s}{p_{c}-1}}w)\rho_{\eta}{\mathrm{d}}y+\Sigma_{1}^{1}(s)+\Sigma_{1}^{2}(s)+\Sigma_{1}^{3}(s).
\end{eqnarray}
where
\begin{eqnarray*}
\Sigma_{1}^{2}(s)&=&\frac{\eta^2(p_{c}+3)}{2}\int_{B}w\partial_{s} w\rho_{\eta} {\mathrm{d}}y,\\ 
\Sigma_{1}^{3}(s)&=&-\eta e^{-\frac{2p_{c}s}{p_{c}-1}}\int_{B}wf(e^{\frac{2s}{p_{c}-1}}w)\rho_{\eta}{\mathrm{d}}y+\frac{\eta(p_{c}+3)}{2}e^{\frac{-2(p_{c}+1)s}{p_{c}-1}}\int_{B} F(e^{\frac{2s}{p_{c}-1}}w)\rho_{\eta} {\mathrm{d}}y.\\
\end{eqnarray*}
We now study each of this last three terms. By the Cauchy-Schwarz inequality, we write for all $\mu\in (0,1)$   
$$\Sigma_{1}^{1}(s)\leq 2\eta(1-\mu)\int_{B}(\partial_{s} w)^2\frac{\rho_{\eta}}{1-|y|^2}{\mathrm{d}}y+\frac{2\eta^3\mu}{1-\mu}\int_{B}w^2\frac{|y|^{2}\rho_{\eta}}{1-|y|^2}{\mathrm{d}}y.$$
We apply the Hardy type inequality to the second term (for the sake of completness, we postpone to Appendix A a short proof) and we choose $\mu= \frac{p_{c}-1}{p_{c}+15}$, we conclude that
\begin{eqnarray}\label{16}
\Sigma_{1}^{1}(s)&\leq &\frac{32\eta}{p_{c}+15}\int_{B}(\partial_{s} w)^2\frac{\rho_{\eta}}{1-|y|^2}{\mathrm{d}}y +\frac{\eta^2N(p_{c}-1)}{8}\int_{B} w^2\rho_{\eta}{\mathrm{d}}y\nonumber \\
&&+\frac{\eta(p_{c}-1)}{8}\int_{B}|\nabla w|^2(1-|y|^2)\rho_{\eta} {\mathrm{d}}y.
\end{eqnarray}
To estimate $\Sigma_{1}^{2}(s)$, we use the Cauchy-Schwarz inequality and we get 
\begin{equation}\label{18}
\Sigma_{1}^{2}(s)\leq \frac{\eta (p_{c}-1)}{8}\int_{B}(\partial_{s} w)^2\rho_{\eta} {\mathrm{d}}y+C\int_{B} w^2\rho_{\eta} {\mathrm{d}}y.
\end{equation}
Using $\eqref{200}$, we obtain for all $s\geq \max (s_{0},1)$
\begin{equation}\label{30}
\Sigma_{1}^{3}(s)\leq \frac{C}{s^a}\int_{B}|w|^{p_{c}+1} \rho_{\eta} {\mathrm{d}}y+Ce^{-\frac{p_{c}+1}{p_{c}-1}s} .
\end{equation}
Finally by using $\eqref{12}$, $\eqref{16}$, $\eqref{18}$ and $\eqref{30}$ we have easily the estimates $\eqref{6}$ and $\eqref{7}$, which ends the proof of Lemma 2.3.
\Box

\bigskip

From Lemmas 2.2 and 2.3, we are in position to prove Proposition 2.1.
\qquad
\\{\it Proof of Proposition 2.1}:
As in \cite{MH1}, we combine Lemmas 2.2 and 2.3, we choose $S_{0}\geq 1$ large enough, so that for all $s\geq \max (s_{0},S_{0})$, $\frac{\eta (p_{c}-1)}{4(p_{c}+1)}-\frac{C}{s^a}\geq 0$ and we use Jensen's inequality to estimate $C \int_{B} w^2\rho_{\eta}dy$, we obtain for all $s\geq \max (s_{0},S_{0})$
\begin{eqnarray*}\label{24}
\frac{d}{ds}(G_{\eta}(w(s),s))\hspace{-0,3cm} &\leq &\hspace{-0,3cm}-e^{\frac{-\eta(p_{c}+3)s}{2}}\Big(\frac{\eta (p_{c}-1)}{p_{c}+15}\int_{B}(\partial_{s} w)^2\frac{\rho_{\eta}}{1-|y|^2}{\mathrm{d}}y+\frac{\eta (p_{c}-1)}{8(p_{c}+1)}\int_{B}|w|^{p_{c}+1}\rho_{\eta} {\mathrm{d}}y \\
&&+\frac{\eta (p_{c}-1)}{8}\int_{B}|\nabla w|^2(1-|y|^2)\rho_{\eta} {\mathrm{d}}y +C-\frac{\eta \Theta (p_{c}+3)}{2}\Big).
\end{eqnarray*}
We now choose $\lambda_{0}=\eta\min(\frac{ p_{c}-1}{p_{c}+15},\frac{ p_{c}-1}{8(p_{c}+1)},\frac{ p_{c}-1}{8})=\frac{ p_{c}-1}{8(p_{c}+1)}$ and $\theta =\theta(\eta)$ large enough, so we have $C-\frac{\eta \theta (p_{c}+3)}{2}\leq 0$, which ends the proof of the first part of Proposition 2.1.
\\To end the proof of the last point of this proposition, we refer the reader to \cite{X}. Note that, our proof strongly relies on the fact that $p_{c}\equiv 1+\frac{4}{N-1}<1+\frac{4}{N-2}$. 
\Box
\subsection{Proof of Theorem 1 }
 We define the following time
$$t_{0}(x_{0})=max(T(x_{0})-e^{-S_{1}},0).$$
Since $\eta \in (0,1)$, according to the Proposition 2.1, we obtain the following corollary which summarizes the principle properties of $H_{\eta}(w(s),s)$.
\begin{cor}{\bf (Estimate on $H_{\eta}(w(s),s)$).}
For all $\eta \in (0,1)$, there exists $t_{0}(x_{0})\in [0,T(x_{0}))$ such that, for all $T_{0}\in (t_{0}(x_{0}),T(x_{0})]$, for all $s \geq -\log (T_{0}-t_{0}(x_{0}))$ and $y\in B$ we have
$$-C\leq H_{\eta}(w(s),s) \leq \Big(\theta +H_{\eta}(w(s_{0}),s_{0}) \Big)e^{\frac{\eta(p_{c}+3)s}{2}},$$
$$\int_{s}^{s+1}\int_{B}(\partial_{s} w(y,\tau))^2\frac{\rho_{\eta}}{1-|y|^2}{\mathrm{d}}y{\mathrm{d}}\tau\leq C\Big(\theta +H_{\eta}(w(s_{0}),s_{0} ) \Big)e^{\frac{\eta(p_{c}+3)s}{2}}, $$
$$ \int_{s}^{s+1}\int_{B_{\frac{1}{2}}}|w(y,\tau)|^{p_{c}+1} {\mathrm{d}}y{\mathrm{d}}\tau +\int_{s}^{s+1}\int_{B_{\frac{1}{2}}}|\nabla w(y,\tau)|^2 {\mathrm{d}}y{\mathrm{d}}\tau\leq C\Big(\theta +H_{\eta}(w(s_{0}),s_{0} ) \Big)e^{\frac{\eta(p_{c}+3)s}{2}}, $$
where $w=w_{x_{0},T_{0}}$ is defined in $\eqref{150}$.

\end{cor}
\begin{rem}
Using the definition $\eqref{150}$ of $w=w_{x_{0},T_{0}}$, we write easily 
$$C\theta +CH_{\eta}(w(s_{0}),s_{0}) \leq K_{0},$$
where $ K_{0}= K_{0}\Big(\eta,T_{0}-t_{0}(x_{0}),\|(U(t_{0}(x_{0})),\partial_{t}U(t_{0}(x_{0})))\|_{H^{1}\times L^{2}(B(x_{0},\frac{T_{0}-t_{0}(x_{0})}{\delta_{0}(x_{0})}))}\Big)$ and $\delta_{0}(x_{0})\in (0,1)$ is defined in $\eqref{17}$.
\end{rem}

From Corollary 2.4, we are in position to prove Theorem 1.\\
{\it Proof of Theorem 1}:
Note that the estimate on the space-time $L^{2}$ norm of $\partial_{s} w$ was already proved in Corollary 2.4. Thus we focus on the space-time $L^{p_{c}+1}$ norm of $w$ and  $L^{2}$ norm of $\nabla w$.\,\,This estimate proved in Corollary 2.4 just for the space-time $L^{p_{c}+1}$ norm of $w$ and  $L^{2}$ norm of $\nabla w$ in $B_{\frac{1}{2}}$.\,\,To extend this estimate from $B_{\frac{1}{2}}$ to $B$ we refer the reader to Merle and Zaag \cite{kl} (unperturbed case) and Hamza and Zaag \cite{MH1} (perturbed case), where they introduce a new covering argument to extend the estimate of any known space $L^{q}$ norm of $w$, $\partial_{s} w$, or $\nabla w$, from $B_{\frac{1}{2}}$ to $B$.
\Box
 
\section{Proof of Theorem 4 }
In this section, we assume that $U$ is a radial blow-up solution of $\eqref{oy}$ with 
\\$g(x,t,\nabla U ,\partial_{t} U )=g(|x|,t,\nabla U.\frac{x}{|x|},\partial_{t} U )$ and $a>1$.
We prove Theorem 4 here, before doing that let us remark that, as we mentionned above the exponential bound is not sufficient to conclude Theorem 6 (unlike Hamza and Zaag \cite{MH1}). According to Theorem 1 we obtain here the polynomially bound on time average of the $H^{1}\times L^{2}((0,1))$ norm of $(w, \partial_{s} w)$ in Theorem 4, this estimate are effective just in the radial case, throughout this bound we can conclude our optimal result written in Theorem 6. 
\\We proceed in three subsections:
\begin{itemize}
\item First, we prove Proposition 3.1, where we show that $L(w(s),s)$ defined in $\eqref{F}$ is a Lyapunov functional for equation $\eqref{0301}$ which is crucial to deduce $(i)$ and $(ii)$ of Theorem 4.
\item Then, we obtain an exponential bound to the time average of the $\int_{0}^{1} |w|^{p_{c}+1}\frac{y^{N-1} \rho_{\eta}}{1-y^2} {\mathrm{d}}y$ when $a> 1$, to do that, we multiply $\eqref{p'}$ by $y\partial_{y} w$ we obtain the Pohozaev identity, which is crucial to prove $(iii)$ of Theorem 4.
\item Finally, the third subsection is devoted to the conclusion of Theorem 4.
\end{itemize}

\subsection{A Lyapunov functional for equation $\eqref{0301}$}
According to the change of variables (\ref{0420}), we can see that the function $w=w_{0}$ satisfies the following equation for all $y\in (0,1)$ and $s\geq -\log(T_{0}):$
\begin{eqnarray}\label{0301}
\partial^2_{s}w &=&\frac{1}{\varphi (y,s)} \partial_{y}(\varphi (y,s)(1-y^2)\partial_{y}w)+\frac{2}{s^{b}}y\partial_{y}w-\frac{2(p_{c}+1)}{(p_{c}-1)^2}w\nonumber \\
&&+|w|^{p_{c}-1}w-\frac{p_{c}+3}{p_{c}-1}\partial_{s}w-2y \partial^2_{y,s}w+e^{\frac{-2p_{c}s}{p_{c}-1}}f(e^{\frac{2s}{p_{c}-1}}w),
\end{eqnarray}
where $\varphi (y,s)$ is defined in $(\ref{phi})$.

\begin{rem}
{\rm  It's worth noticing to recall that the weight $\rho_{s^{-b}}$ which is defined in $\eqref{phi}$ depends on time, it is not the case in this series of papers \cite{MH}, \cite{MH1}, \cite{MH2}, \cite{X}, \cite{fh6}, \cite{fh7}, \cite{fh3}, \cite{fh4}, \cite{kl}, \cite{fh} and \cite{fh8}, we expect that the derivations in time is problematic, in fact, we note after observation as we said above, that there are new terms appearing compared to the previous works which was already controlled.}
\end{rem}
The equation $\eqref{0301}$ will be studied in the following space $\H_{rad}$
$$\H_{rad}=\{q=(q_{1},q_{2})\in H^{1}_{loc,u}\times L^{2}_{loc,u}((0,1))| \int_{0}^{1} \Big(q_{2}^2+q'^2_{1}(1-y^2)+q_{1}^2\Big)\rho {\mathrm{d}}y< +\infty\}. $$
To control the norm of $(w(s),\partial_{s} w(s) ) \in \H_{rad}$, we first introduce the following functionals:
\begin{equation}\label{F}
 L(w(s),s)=\exp\Big(\frac{p_{c}+3}{2(b-1)s^{b-1}}\Big) K(w(s),s)+\frac{\sigma}{s^{b-1}}\,\,\,\,{\rm with }\,\,\,b\in (1,a),\\
\end{equation}
and
\begin{eqnarray}\label{FF}
K(w(s),s)&=&E(w(s),s)+J(w(s),s),\nonumber\\
E(w(s),s) \hspace{-0,3cm}&=& \hspace{-0,4cm}\int_{0}^{1}\Big(\frac{1}{2}(\partial_{s}w)^2+\frac{1}{2}(\partial_{y}w)^2(1-y^2)+\frac{p_{c}+1}{(p_{c}-1)^2}w^2-\frac{1}{p_{c}+1}|w|^{p_{c}+1}\Big)\varphi (y,s) {\mathrm{d}}y\nonumber\\
&&-e^{\frac{-2(p_{c}+1)s}{p_{c}-1}} \int_{0}^{1} F(e^{\frac{2s}{p_{c}-1}}w)\varphi (y,s) {\mathrm{d}}y,\\
J(w(s),s)&=&-\frac{1}{s^{b}}\int_{0}^{1}  w\partial_{s}w \varphi (y,s){\mathrm{d}}y,\nonumber
\end{eqnarray}
and $\sigma$ is a constant will be determined later. We are going to prove that $L(w(s),s)$ is the Lyapunov functional of equation $\eqref{0301}$ in the following proposition then we can deduce directly $(i)$ and $(ii)$ of Theorem 4.
\begin{prop}
For all $b\in (1,a)$, there exists $S_{4}\geq 1$ and $\lambda_{1}>0$ such that $L(w(s),s)$ defined in $\eqref{F}$ satisfies for all $s\geq \max(S_{4},s_{0})$
 \begin{eqnarray}\label{new}
L(w(s+1),s+1)-L(w(s),s)&\leq &-\frac{\lambda_{1}}{s^{b}}\int_{s}^{s+1}\int_{0}^{1}|w(y,\tau)|^{p_{c}+1}\varphi (y,\tau) {\mathrm{d}}y{\mathrm{d}}\tau \nonumber\\
&&-\frac{\lambda_{1}}{s^{b}}\int_{s}^{s+1}\int_{0}^{1}(\partial_{s}w(y,\tau))^2\frac{y^{2}\varphi (y,\tau)}{1-y^2}{\mathrm{d}}y{\mathrm{d}}\tau \\
&&-\frac{\lambda_{1}}{s^{b}}\int_{s}^{s+1}\int_{0}^{1}(\partial_{y}w(y,\tau))^2(1-y^2)\varphi (y,\tau) {\mathrm{d}}y{\mathrm{d}}\tau ,\nonumber
\end{eqnarray}
where $w= w_{0,T_{0}}$ defined in $\eqref{0420}$. 
\\Moreover, there exists $S_{5}\geq S_{4}$ such that for all $s\geq max(S_{5},s_{0})$ $L(w(s),s)\geq 0$.
\end{prop}

The existence of this Lyapunov functional (and a blow-up criterion for equation $\eqref{0301}$ based on $L(w(s),s)$) are a crucial step in the derivation of Theorem 4. Indeed with the functional $H_{0}(w(s),s)$ (defined below in $\eqref{kdv1}$) and some more works, we are able to adapt the analysis performed in \cite{kl} for equation $\eqref{1.6}$ to get Theorem 6.
We aim at proving that the functional $L(w(s),s)$ defined in $\eqref{F}$ is a Lyapunov functional for equation $\eqref{0301}$, provided that $s$ is large enough. 
 We give here the following result:
\begin{lem}
For all $b\in (1,a)$ and $\varepsilon_{1}\in (0,\frac{1}{2})$, there exists $S_{2}\geq S_{1}$ such that we have for all $s \geq \max(S_{2},s_{0})$
 \begin{eqnarray*}
\frac{d}{ds}(E(w(s),s))&\leq &-\frac{2}{s^{b}}\int_{0}^{1}(\partial_{s}w)^2\frac{y^{2}\varphi (y,s)}{1-y^2}{\mathrm{d}}y+\frac{2}{s^{b}}\int_{0}^{1}\partial_{s}w\partial_{y}wy\varphi (y,s) {\mathrm{d}}y\\
&&+\frac{1}{4(p_{c}+1)s^{b+1}}\int_{0}^{1}|w|^{p_{c}+1}\log(1-y^2)\varphi (y,s){\mathrm{d}}y+\Sigma_{2}(s),
\end{eqnarray*}
where $\Sigma_{2}(s)$ satisfies
    \begin{eqnarray*}\label{a4}
\Sigma_{2}(s)&\leq & \frac{p_{c}-1}{16s^{b}} \int_{0}^{1}(\partial_{y}w)^2\varphi (y,s)(1-y^2){\mathrm{d}}y+\frac{p_{c}-1}{16s^{b}}\int_{0}^{1}(\partial_{s}w)^2\varphi (y,s){\mathrm{d}}y\\
&&+\frac{C}{s^a}\int_{0}^{1}|w|^{p_{c}+1} \varphi (y,s) {\mathrm{d}}y + C e^{-\frac{p_{c}-1}{32b} s}\int_{0}^{1}\frac{(\partial_{s}w)^2}{1-y^2}y^{N-1}(1-y^2)^{\varepsilon_{1}}{\mathrm{d}}y\\
&&+C e^{-\frac{p_{c}-1}{16b} s}\int_{0}^{1}(\partial_{y}w)^2y^{N-1}{\mathrm{d}}y+\frac{C}{s^{b}}.
\end{eqnarray*}
\end{lem}
\bigskip

{\it Proof}: By virtue of identity $\eqref{0301}$, it can be seen that
$$\frac{d}{ds}(E(w(s),s))=-\frac{2}{s^{b}}\int_{0}^{1}(\partial_{s}w)^2\frac{y^{2}\varphi (y,s)}{1-y^2}{\mathrm{d}}y+\frac{2}{s^{b}}\int_{0}^{1}\partial_{s}w\partial_{y}wy\varphi (y,s) {\mathrm{d}}y$$
\begin{equation}\label{XXX} 
+\Sigma_{2}^{1}(s)+\Sigma_{2}^{2}(s)+\Sigma_{2}^{3}(s)+\Sigma_{2}^{4}(s)+\Sigma_{2}^{5}(s)+\Sigma_{2}^{6}(s),\\
 \end{equation}
with 
\begin{eqnarray*}
\Sigma_{2}^{1}(s)&=&\frac{b}{s^{b+1}(p_{c}+1)}\int_{0}^{1}|w|^{p_{c}+1}\log(1-y^2)\varphi (y,s){\mathrm{d}}y,\\
\Sigma_{2}^{2}(s)&=&\hspace{-0,3cm}\frac{2(p_{c}+1)}{p_{c}-1}e^{-\frac{2(p_{c}+1)s}{p_{c}-1}}\int_{0}^{1} F(e^{\frac{2s}{p_{c}-1}}w) \varphi (y,s){\mathrm{d}}y-\frac{2e^{-\frac{2p_{c}s}{p_{c}-1}}}{p_{c}-1}\int_{0}^{1}wf(e^{\frac{2s}{p_{c}-1}}w)\varphi (y,s){\mathrm{d}}y,\\
\Sigma_{2}^{3}(s)&=&-\frac{b}{2s^{b+1}}\int_{0}^{1}(\partial_{s}w)^2\log(1-y^2)\varphi (y,s){\mathrm{d}}y,\\
\Sigma_{2}^{4}(s)&=&-\frac{b}{2s^{b+1}}\int_{0}^{1}(\partial_{y}w)^2(1-y^2)\log(1-y^2)\varphi (y,s){\mathrm{d}}y,\\
\Sigma_{2}^{5}(s)&=&-\frac{b(p_{c}+1)}{s^{b+1}(p_{c}-1)^2}\int_{0}^{1}w^2\log(1-y^2)\varphi (y,s){\mathrm{d}}y,\\
\Sigma_{2}^{6}(s)&=&\frac{be^{-\frac{2(p_{c}+1)s}{p_{c}-1}}}{s^{b+1}}\int_{0}^{1}F(e^{\frac{2s}{p_{c}-1}}w)\log(1-y^2)\varphi (y,s){\mathrm{d}}y.
\end{eqnarray*}
A similar study as $\Sigma_{0}(s)$ defined in $\eqref{ao}$, just we need to replace $\rho_{\eta}=(1-|y|^2)^{\eta}$ by $\rho_{s^{-b}}=(1-y^2)^{s^{-b}}$, gives rise to the following inequality
 \begin{equation}\label{pi}
  \Sigma_{2}^{2}(s) \leq \frac{C}{s^a}\int_{0}^{1}|w|^{p_{c}+1}\varphi (y,s) {\mathrm{d}}y+Ce^{-\frac {p_{c}+1 }{p_{c}-1}s} .\\
\end{equation}
We are going now to estimate $\Sigma_{2}^{3}(s),$ we divide the interval $(0,1)$ into two parts
 \begin{equation}\label{keep}
 B_{1}(s)=\{y \in (0,1)\,|\, 1-y^2\leq  e^{-\frac{p_{c}-1}{8b} s}\},\,B_{2}(s)=\{y \in (0,1)\,|\, 1-y^2>  e^{-\frac{p_{c}-1}{8b} s}\}.\\
\end{equation}
We see easily that
$$\Sigma_{2}^{3}(s)=\chi_{1}(s)+\chi_{2}(s),$$
with 
\begin{eqnarray*}
\chi_{1}(s)&=&-\frac{b}{2s^{b+1}}\int_{B_{1}(s)}(\partial_{s}w)^2\log(1-y^2)\varphi (y,s){\mathrm{d}}y,\\
\chi_{2}(s)&=&-\frac{b}{2s^{b+1}}\int_{B_{2}(s)}(\partial_{s}w)^2\log(1-y^2)\varphi (y,s){\mathrm{d}}y.\\
\end{eqnarray*}
By combining the fact that, the function $y\mapsto (1-y^2)^{\frac{1}{4}}\log(1-y^2)$ is bounded in $(0,1)$ and if $y \in B_{1}(s)$, we have $(1-y^2)^{\frac{1}{4}}\leq  e^{-\frac{p_{c}-1}{32b} s}$, we obtain
$$\chi_{1}(s) \leq \frac{C e^{-\frac{p_{c}-1}{32b} s}}{s^{b+1}}\int_{B_{1}(s)}\frac{(\partial_{s}w)^2}{(1-y^2)^{\frac{1}{2}}}y^{N-1}{\mathrm{d}}y \leq  C e^{-\frac{p_{c}-1}{32b} s}\int_{0}^{1}\frac{(\partial_{s}w)^2}{1-y^2}y^{N-1}(1-y^2)^{\varepsilon_{1}}{\mathrm{d}}y$$
for all $\varepsilon_{1} \in (0,\frac{1}{2})$.
\\If $y \in B_{2}(s)$, we can see that $-\log(1-y^2)\leq \frac{p_{c}-1}{8b}s$, so
$$\chi_{2}(s) \leq \frac{p_{c}-1}{16s^{b}}\int_{B_{2}(s)}(\partial_{s}w)^2\varphi (y,s){\mathrm{d}}y\leq \frac{p_{c}-1}{16s^{b}}\int_{0}^{1}(\partial_{s}w)^2\varphi (y,s){\mathrm{d}}y .$$
We can deduce that, for all $\varepsilon_{1} \in (0,\frac{1}{2})$
\begin{equation}\label{X}
\Sigma_{2}^{3}(s) \leq \frac{p_{c}-1}{16s^{b}}\int_{0}^{1}(\partial_{s}w)^2\varphi (y,s){\mathrm{d}}y+ C e^{-\frac{p_{c}-1}{32b} s}\int_{0}^{1}\frac{(\partial_{s}w)^2}{1-y^2}y^{N-1}(1-y^2)^{\varepsilon_{1}}{\mathrm{d}}y.\\
\end{equation}
Applying the same lines of reasoning as in the treatment of the last term to estimate $\Sigma_{2}^{3}(s),$ (we keep the same partition of the interval $(0,1)$ as $\eqref{keep}$), then we write
$$\Sigma_{2}^{4}(s)=\chi_{3}(s)+\chi_{4}(s),$$
with 
\begin{eqnarray*}
\chi_{3}(s)&=&-\frac{b}{2s^{b+1}}\int_{B_{1}(s)}(\partial_{y}w)^2(1-y^2)\log(1-y^2)\varphi (y,s){\mathrm{d}}y,\\
\chi_{4}(s)&=&-\frac{b}{2s^{b+1}}\int_{B_{2}(s)}(\partial_{y}w)^2(1-y^2)\log(1-y^2)\varphi (y,s){\mathrm{d}}y.\\
\end{eqnarray*}
By combining the fact that, the function $y\mapsto (1-y^2)^{\frac{1}{2}}\log(1-y^2)$ is bounded in $(0,1)$ and if $y \in B_{1}(s)$, we have $(1-y^2)^{\frac{1}{2}}\leq  e^{-\frac{p_{c}-1}{16b} s}$, yields to
 $$\chi_{3}(s) \leq C e^{-\frac{p_{c}-1}{16b} s}\int_{B_{1}(s)}(\partial_{y}w)^2y^{N-1}{\mathrm{d}}y\leq C e^{-\frac{p_{c}-1}{16b} s}\int_{0}^{1}(\partial_{y}w)^2y^{N-1}{\mathrm{d}}y,$$
If $y \in B_{2}(s)$, we know that $-\log(1-y^2)\leq \frac{p_{c}-1}{8b}s$, so
$$\chi_{4}(s) \leq \frac{p_{c}-1}{16s^{b}}\int_{B_{2}(s)}(\partial_{y}w)^2(1-y^2)\varphi (y,s){\mathrm{d}}y\leq \frac{p_{c}-1}{16s^{b}}\int_{0}^{1}(\partial_{y}w)^2(1-y^2)\varphi (y,s){\mathrm{d}}y, $$
which ensure that 
\begin{equation}\label{XXXX}
\Sigma_{2}^{4}(s)\leq \frac{p_{c}-1}{16s^{b}}\int_{0}^{1}(\partial_{y}w)^2(1-y^2)\varphi (y,s){\mathrm{d}}y+C e^{-\frac{p_{c}-1}{16b} s}\int_{0}^{1}(\partial_{y}w)^2y^{N-1}{\mathrm{d}}y.\\
\end{equation}
To estimate $\Sigma_{2}^{4}(s)$, we start by recalling the following Young inequality:
\begin{equation}\label{502}
 w^2 \leq \frac{1}{\varepsilon}+ \varepsilon|w|^{p_{c}+1}\,\,\,\,\,\,\,\,\forall\,\,\,\,\varepsilon> 0,\\
\end{equation} 
We choose $\varepsilon =\frac{(p_{c}-1)^2}{2(p_{c}+1)^2}$, then we multiply $\eqref{502}$ by $-\frac{(p_{c}+1)b}{(p_{c}-1)^2s^{b+1}}\log(1-y^2)\varphi (y,s)$ and integrate over $(0,1)$, the simple fact that $-\int_{0}^{1}\log(1-y^2)\varphi (y,s){\mathrm{d}}y \leq C$ gives rise to the following result:
\begin{equation}\label{503}
\Sigma_{2}^{5}(s) \leq 
-\frac{b}{2(p_{c}+1)s^{b+1}}\int_{0}^{1}|w|^{p_{c}+1}\log(1-y^2)\varphi (y,s){\mathrm{d}}y +\frac{C}{s^{b+1}}.\\
\end{equation}
Performing to inequality $\eqref{104}$, we find that
\begin{equation}\label{506}
e^{-\frac{2(p_{c}+1)s}{p_{c}-1}}|F(e^{\frac{2s}{p_{c}-1}}w)|\leq C \frac{|w|^{p_{c}+1}}{(\log(2+e^{\frac{4s}{p_{c}-1}}w^2))^a}+Ce^{\frac{-2(p_{c}+1)s}{p_{c}-1}}.\\
\end{equation}
Multiplying $\eqref{506}$ by $-\frac{b}{s^{b+1}}\log(1-y^2)\varphi (y,s)$ and integrate over $(0,1)$, 
we infer
 $$\Sigma_{2}^{6}(s) \leq -\frac{C}{s^{b+1}}\int_{0}^{1}\frac{|w|^{p_{c}+1}}{(\log(2+e^{\frac{4s}{p_{c}-1}}w^2))^a}\log(1-y^2)\varphi (y,s) {\mathrm{d}}y+\frac{C}{s^{b+1}},$$
we deduce in view of $\eqref{200}$
\begin{equation}\label{509}
\Sigma_{2}^{6}(s)\leq \frac{C}{s^{b+1}}-\frac{C}{s^{a+b+1}}\int_{0}^{1}|w|^{p_{c}+1} \log(1-y^2)\varphi (y,s){\mathrm{d}}y.\\
\end{equation}
Now, we combine $\eqref{503}$ and $\eqref{509}$ to write 
\begin{eqnarray}\label{550}
\Sigma_{2}^{1}(s)+\Sigma_{2}^{5}(s)+\Sigma_{2}^{6}(s)&\leq &-\Big(\frac{C}{s^{a}}-\frac{b}{2(p_{c}+1)}\Big)\frac{1}{s^{b+1}}\int_{0}^{1}|w|^{p_{c}+1} \log(1-y^2)\varphi (y,s){\mathrm{d}}y\nonumber \\
&&+\frac{C}{s^{b+1}},
\end{eqnarray}
from $\eqref{550}$ we can see that there exists $S_{2}\geq S_{1}$ such that we have for all $s \geq \max(S_{2},s_{0})$:
\begin{equation}\label{551}
\Sigma_{2}^{1}(s)+\Sigma_{2}^{5}(s)+\Sigma_{2}^{6}(s)\leq \frac{b}{4(p_{c}+1)s^{b+1}}\int_{0}^{1}|w|^{p_{c}+1} \log(1-y^2)\varphi (y,s){\mathrm{d}}y+\frac{C}{s^{b+1}}.\\
\end{equation}
The result derives immediately from $\eqref{XXX}$, $\eqref{pi}$, $\eqref{X}$, $\eqref{XXXX}$ and $\eqref{551}$, which ends the proof of Lemma 3.2.
\Box

We are going to prove the following estimate to the functional $J(w(s),s)$. 
 \begin{lem}
For all $b\in (1,a)$ and $\varepsilon_{1}\in (0,\frac{1}{2})$, there exists $S_{3}\geq S_{2}$ such that we have for all $s \geq \max(S_{3},s_{0})$
\begin{eqnarray*}
\frac{d}{ds}(J(w(s),s))&\leq & \frac{32}{(p_{c}+15)s^{b}}\int_{0}^{1}(\partial_{s}w)^2\frac{y^{2}\varphi (y,s)}{1-y^2}{\mathrm{d}}y+\frac{p_{c}+3}{2s^{b}}K(w(s),s)\\
&-&\frac{2}{s^{b}}\int_{0}^{1} \partial_{y}w\partial_{s}w y\varphi (y,s){\mathrm{d}}y-\frac{p_{c}-1}{8s^{b}}\int_{0}^{1}(\partial_{y}w)^2(1-y^2)\varphi (y,s){\mathrm{d}}y\\
&-&\hspace{-0,4cm}\frac{p_{c}+7}{16s^{b}}\int_{0}^{1}\hspace{-0,2cm}(\partial_{s}w)^2\varphi (y,s) {\mathrm{d}}y-\frac{p_{c}-1}{8(p_{c}+1)s^{b}}\int_{0}^{1}\hspace{-0,2cm}|w|^{p_{c}+1}\varphi (y,s){\mathrm{d}}y+\hspace{-0,2cm}\Sigma_{3}(s),
\end{eqnarray*}
where $\Sigma_{3}(s)$ satisfies
$$\Sigma_{3}(s)\leq -\frac{b}{s^{2b+1}}\int_{0}^{1}\hspace{-0,2cm}|w|^{p_{c}+1}\log(1-y^2)\varphi (y,s){\mathrm{d}}y+C e^{-\frac{p_{c}-1}{32b} s}\int_{0}^{1}\hspace{-0,2cm}\frac{(\partial_{s}w)^2}{1-y^2}y^{N-1}(1-y^2)^{\varepsilon_{1}}{\mathrm{d}}y+\frac{C}{s^{b}}.$$
\end{lem} 

{\it Proof}: Note that $J(w(s),s)$ is a differentiable function according to equation $\eqref{0301}$ we get for all $s\geq \max(s_{0},1)$ 
\begin{eqnarray*}
\frac{d}{ds}(J(w(s),s))&\leq &-\frac{1}{s^{b}}\int_{0}^{1} (\partial_{s}w)^2 \varphi (y,s) {\mathrm{d}}y+\frac{1}{s^{b}}\int_{0}^{1}(\partial_{y}w)^2(1-y^2)\varphi (y,s){\mathrm{d}}y\nonumber \\
&&+\frac{2(p_{c}+1)}{(p_{c}-1)^2s^{b} }\int_{0}^{1} w^2 \varphi (y,s) {\mathrm{d}}y-\frac{1}{s^{b}}\int_{0}^{1}|w|^{p_{c}+1}\varphi (y,s) {\mathrm{d}}y\nonumber \\
 && +\Big(\frac{b}{s}-2N+\hspace{-0,2cm}\frac{p_{c}+3}{p_{c}-1}\Big)\frac{1}{s^{b}} \hspace{-0,2cm}\int_{0}^{1} \hspace{-0,3cm}w\partial_{s}w \varphi (y,s){\mathrm{d}}y-\frac{2}{s^{b}}\hspace{-0,2cm}\int_{0}^{1}\hspace{-0,3cm} \partial_{y}w\partial_{s}w y\varphi (y,s){\mathrm{d}}y\nonumber \\
 &&+\frac{4}{s^{2b}}\int_{0}^{1}w \partial_{s}w\frac{y^{2}\varphi (y,s)}{1-y^2} {\mathrm{d}}y -\frac{e^{\frac{-2p_{c}s}{p_{c}-1}}}{s^{b}}\int_{0}^{1} wf(e^{\frac{2s}{p_{c}-1}}w)\varphi (y,s) {\mathrm{d}}y \nonumber \\
&&-\frac{ 2}{s^{2b}}\int_{0}^{1}w\partial_{y}wy\varphi (y,s){\mathrm{d}}y+\frac{b}{s^{2b+1}}\int_{0}^{1}w\partial_{s}w\log(1-y^2)\varphi (y,s){\mathrm{d}}y.
\end{eqnarray*} 
According to the expression of $K(w(s),s)$ in $\eqref{FF}$, with some straighforward computation we show the following inequality
 \begin{eqnarray}\label{801}
\frac{d}{ds}(J(w(s),s))&\leq &-\frac{p_{c}+7}{4s^{b}}\int_{0}^{1}(\partial_{s}w)^2 \varphi (y,s){\mathrm{d}}y +\frac{p_{c}+3}{2s^{b}}K(w(s),s)  \\
&&-\frac{p_{c}-1}{4s^{b}}\int_{0}^{1}(\partial_{y}w)^2(1-y^2)\varphi (y,s) {\mathrm{d}}y-\frac{p_{c}-1}{2(p_{c}+1)s^{b}}\int_{0}^{1}|w|^{p_{c}+1}\varphi (y,s) {\mathrm{d}}y\nonumber \\
&&-\frac{2}{s^{b}}\int_{0}^{1} \partial_{y}w\partial_{s}w y\varphi (y,s){\mathrm{d}}y+\Sigma_{3}^{1}(s)+\Sigma_{3}^{2}(s)+\Sigma_{3}^{3}(s)+\Sigma_{3}^{4}(s)+\Sigma_{3}^{5}(s),\nonumber 
  \end{eqnarray}  
 such that
  \begin{eqnarray*}
 \Sigma_{3}^{1}(s)&=&\Big(\frac{N}{s^{b}}-\frac{p_{c}+1}{2(p_{c}-1)}\Big)\frac{1}{s^{b}}\int_{0}^{1}w^2\varphi (y,s){\mathrm{d}}y,\\
 \Sigma_{3}^{2}(s)&=&\frac{b}{s^{2b+1}}\int_{0}^{1}w\partial_{s}w\log(1-y^2)\varphi (y,s){\mathrm{d}}y,\\
  \Sigma_{3}^{3}(s)&=&\Big(\frac{b}{s}-N+\frac{p_{c}+3}{2s^{b}}\Big)\frac{1}{s^{b}}  \int_{0}^{1} w\partial_{s}w \varphi (y,s) {\mathrm{d}}y,\\
  \Sigma_{3}^{4}(s)&=&\frac{4}{s^{2b}}\int_{0}^{1}w \partial_{s}w\frac{y^{2}\varphi (y,s)}{1-y^2} {\mathrm{d}}y-\frac{2}{s^{3b}}\int_{0}^{1}w^2\frac{y^{2}\varphi (y,s)}{1-y^2}{\mathrm{d}}y,\\
  \Sigma_{3}^{5}(s)&=&-\frac{e^{\frac{-2p_{c}s}{p_{c}-1}}}{s^{b}}\int_{0}^{1} wf(e^{\frac{2s}{p_{c}-1}}w)\varphi (y,s) {\mathrm{d}}y+\frac{p_{c}+3}{2s^{b}}e^{\frac{-2(p_{c}+1)s}{p_{c}-1}} \int_{0}^{1} F(e^{\frac{2s}{p_{c}-1}}w)\varphi (y,s) {\mathrm{d}}y.
  \end{eqnarray*}
We are going now to estimate each of these last five terms, the Cauchy-Schwarz inequality implies that
 \begin{equation}\label{900}
 \Sigma_{3}^{2}(s)\leq -\frac{b}{s^{2b+1}}\int_{0}^{1}w^2\log(1-y^2)\varphi (y,s){\mathrm{d}}y
   +\Sigma_{2}^{2}(s).\\
  \end{equation}
Combining the Young inequality with $\eqref{X}$ and $\eqref{900}$, we obtain for all $\varepsilon_{1}\in (0,\frac{1}{2})$
\begin{eqnarray}\label{0802}
   \Sigma_{3}^{2}(s)&\leq &\frac{C}{s^{2b}}\int_{0}^{1}(\partial_{s}w)^2\varphi (y,s){\mathrm{d}}y-\frac{b}{s^{2b+1}}\int_{0}^{1}|w|^{p_{c}+1}\log(1-y^2)\varphi (y,s){\mathrm{d}}y\nonumber \\
   &&+C e^{-\frac{p_{c}-1}{32b} s}\int_{0}^{1}\frac{(\partial_{s}w)^2}{1-y^2}y^{N-1}(1-y^2)^{\varepsilon_{1}}{\mathrm{d}}y+\frac{C}{s^{2b+1}}.
  \end{eqnarray}
Using the fact that for all $s \geq \max(s_{0},1)$ $\Big|\frac{b}{s}-N+\frac{p_{c}+3}{2s^{b}}\Big|\leq C,$ we get by virtue of the Cauchy-Schwarz inequality    
  \begin{eqnarray}\label{R}
  \Sigma_{3}^{3}(s)\leq \frac{p_{c}+7}{8s^{b}}\int_{0}^{1}(\partial_{s}w)^2\varphi (y,s){\mathrm{d}}y + \frac{\alpha (b,N)} {s^{b}}\int_{0}^{1}w^2\varphi (y,s){\mathrm{d}}y.
  \end{eqnarray}
We choose $\varepsilon =\frac{p_{c}-1}{4\alpha (b,N)(p_{c}+1)}$ in $\eqref{502}$ and according to  $\eqref{R}$, we can deduce that for all $s \geq \max(s_{0},1)$
\begin{equation}\label{802}
   \Sigma_{3}^{3}(s)\leq \frac{p_{c}+7}{8s^{b}}\int_{0}^{1}(\partial_{s}w)^2\varphi (y,s){\mathrm{d}}y + \frac{p_{c}-1}{4(p_{c}+1)s^{b}}\int_{0}^{1}|w|^{p_{c}+1}\varphi (y,s) {\mathrm{d}}y+\frac{C}{s^{b}}.\\
  \end{equation}
By the Cauchy-Schwarz inequality, we write for all $\mu\in (0,1)$ 
\begin{equation}\label{N1}
\Sigma_{3}^{4}(s)\leq \frac{2}{s^{b}}(1-\mu)\int_{0}^{1}(\partial_{s}w)^2\frac{y^{2}\varphi (y,s)}{1-y^2}{\mathrm{d}}y+\frac{2\mu}{(1-\mu)s^{3b}}\int_{0}^{1}w^2\frac{y^{2}\varphi (y,s)}{1-y^2}{\mathrm{d}}y.\\
\end{equation}
By exploiting the Hardy-Sobolev inequality $\eqref{A}$, we get
\begin{equation}\label{N2}
\int_{0}^{1}w^2\frac{y^{2}\varphi (y,s)}{1-y^2}{\mathrm{d}}y\leq s^{2b}\int_{0}^{1}(\partial_{y}w)^2\varphi (y,s)(1-y^2){\mathrm{d}}y+Ns^{b}\int_{0}^{1}w^2\varphi (y,s){\mathrm{d}}y,\\
\end{equation}
from $\eqref{N1}$, $\eqref{N2}$ and if we choose $\mu= \frac{p_{c}-1}{p_{c}+15}$, we conclude that
\begin{eqnarray}\label{920}
\Sigma_{3}^{4}(s)&\leq &\frac{32}{(p_{c}+15)s^{b}}\int_{0}^{1}(\partial_{s} w)^2\frac{y^{2}\varphi (y,s)}{1-y^2}{\mathrm{d}}y +\frac{N(p_{c}-1)}{8s^{2b}}\int_{0}^{1} w^2\varphi (y,s) {\mathrm{d}}y\nonumber \\
&&+\frac{p_{c}-1}{8s^{b}}\int_{0}^{1}(\partial_{y}w)^2\varphi (y,s)(1-y^2){\mathrm{d}}y.
\end{eqnarray}
By $\eqref{920}$, we can see easly that:
\begin{eqnarray}\label{555}
\Sigma_{3}^{1}(s)+\Sigma_{3}^{4}(s)&\leq &\hspace{-0,3cm}\frac{32}{(p_{c}+15)s^{b}}\int_{0}^{1}(\partial_{s} w)^2\frac{y^{2}\varphi (y,s)}{1-y^2}{\mathrm{d}}y +\frac{p_{c}-1}{8s^{b}}\hspace{-0,2cm}\int_{0}^{1}\hspace{-0,3cm}(\partial_{y}w)^2\varphi (y,s)(1-y^2){\mathrm{d}}y\nonumber \\
&&+\Big(\frac{C}{s^{b}}-\frac{p_{c}+1}{2(p_{c}-1)}\Big)\frac{1}{s^{b}}\int_{0}^{1} w^2\varphi (y,s) {\mathrm{d}}y.
\end{eqnarray}
Finally, we estimate $\Sigma_{3}^{5}(s)$ by using inequality $\eqref{pi}$
 \begin{equation}\label{808}
\Sigma_{3}^{5}(s)\leq \frac{C}{s^{a+b}} \int_{0}^{1}|w|^{p_{c}+1}\varphi (y,s) {\mathrm{d}}y+\frac{C}{s^{b}}.
\end{equation}
Combining $\eqref{801}$, $\eqref{0802}$, $\eqref{802}$, $\eqref{555}$ and $\eqref{808}$, we obtain 
\begin{eqnarray*}
\frac{d}{ds}(J(w(s),s))&\leq & \frac{p_{c}+3}{2s^{b}}K(w(s),s)+\frac{32}{(p_{c}+15)s^{b}}\int_{0}^{1}(\partial_{s}w)^2\frac{y^{2}\varphi (y,s)}{1-y^2}{\mathrm{d}}y\\
&-&\frac{2}{s^{b}}\int_{0}^{1} \partial_{y}w\partial_{s}w y\varphi (y,s){\mathrm{d}}y-\frac{p_{c}-1}{8s^{b}}\int_{0}^{1}(\partial_{y}w)^2(1-y^2)\varphi (y,s){\mathrm{d}}y\\
&-&\frac{b}{s^{2b+1}}\int_{0}^{1}|w|^{p_{c}+1}\log(1-y^2)\varphi (y,s){\mathrm{d}}y \\
&+&C e^{-\frac{p_{c}-1}{32b} s}\int_{0}^{1}\frac{(\partial_{s}w)^2}{1-y^2}y^{N-1}(1-y^2)^{\varepsilon_{1}}{\mathrm{d}}y\\
&+&\Big(\frac{C}{s^{b}}-\frac{p_{c}+7}{8}\Big)\frac{1}{s^{b}}\int_{0}^{1}(\partial_{s}w)^2\varphi (y,s) {\mathrm{d}}y\\
&+&\Big(\frac{C}{s^{b}}-\frac{p_{c}+1}{2(p_{c}-1)}\Big)\frac{1}{s^{b}}\int_{0}^{1}w^2\varphi (y,s){\mathrm{d}}y\\
&+&\Big(\frac{C}{s^{a}}-\frac{p_{c}-1}{4(p_{c}+1)}\Big)\frac{1}{s^{b}}\int_{0}^{1}|w|^{p_{c}+1}\varphi (y,s){\mathrm{d}}y +\frac{C}{s^{b}}. 
\end{eqnarray*}
We choose $S_{3}\geq S_{2}$ such that we have $\forall\,\,s\geq \max(S_{3},s_{0})$
 $$\frac{C}{s^{b}}-\frac{p_{c}+7}{16}\leq 0,\,\,\,\frac{C}{s^{b}}-\frac{p_{c}+1}{2(p_{c}-1)}\leq 0,\,\,\,\frac{C}{s^{a}}-\frac{p_{c}-1}{8(p_{c}+1)}\leq 0,$$
 which ends the proof of Lemma 3.3.
\Box

\medskip

Lemmas 3.2 and 3.3 allows to prove Proposition 3.1.\\
{\it Proof of Proposition 3.1}:
 Combining Lemmas 3.2 and 3.3 we can deduce that for all $s\geq \max(S_{3},s_{0})$ and $\varepsilon_{1}\in (0,\frac{1}{2})$
  \begin{eqnarray*}
\frac{d}{ds}(K(w(s),s))&\leq & \frac{p_{c}+3}{2s^{b}}K(w(s),s)-\frac{2(p_{c}-1)}{(p_{c}+15)s^{b}}\int_{0}^{1}(\partial_{s}w)^2\frac{y^{2}\varphi (y,s)}{1-y^2}{\mathrm{d}}y\\
&-&\frac{p_{c}-1}{16s^{b}}\int_{0}^{1}(\partial_{y}w)^2(1-y^2)\varphi (y,s) {\mathrm{d}}y\\
&+&\hspace{-0,3cm} C e^{-\frac{p_{c}-1}{32b} s}\int_{0}^{1}\frac{(\partial_{s}w)^2}{1-y^2}y^{N-1}(1-y^2)^{\varepsilon_{1}}{\mathrm{d}}y+C e^{-\frac{p_{c}-1}{16b} s}\hspace{-0,2cm}\int_{0}^{1}\hspace{-0,3cm}(\partial_{y}w)^2y^{N-1}{\mathrm{d}}y\\
&+&\Big(\frac{C}{s^{a-b}}-\frac{p_{c}-1}{8(p_{c}+1)}\Big)\frac{1}{s^{b}}\int_{0}^{1}|w|^{p_{c}+1}\varphi (y,s) {\mathrm{d}}y \\
&-&\Big(\frac{b}{s^{b}}-\frac{1}{4(p_{c}+1)}\Big)\frac{1}{s^{b+1}}\int_{0}^{1}|w|^{p_{c}+1}\log(1-y^2)\varphi (y,s){\mathrm{d}}y+\frac{C}{s^{b}}.
\end{eqnarray*}
 If we choose $S_{4}\geq S_{3}$ large enough so that $\forall\,\,s\geq \max(S_{4},s_{0})$ we have
 $$\frac{C}{s^{a-b}}-\frac{p_{c}-1}{16(p_{c}+1)}\leq 0,\,\,\,\frac{b}{s^{b}}-\frac{1}{4(p_{c}+1)}\leq 0,$$
 this gives rise to
\begin{eqnarray}\label{ne}
\frac{d}{ds}(K(w(s),s))&\leq & \frac{p_{c}+3}{2s^{b}}K(w(s),s)-\frac{2(p_{c}-1)}{(p_{c}+15)s^{b}}\int_{0}^{1}(\partial_{s}w)^2\frac{y^{2}\varphi (y,s)}{1-y^2}{\mathrm{d}}y\\
&&-\frac{p_{c}-1}{16s^{b}}\int_{0}^{1}(\partial_{y}w)^2(1-y^2)\varphi (y,s) {\mathrm{d}}y-\frac{p_{c}-1}{16(p_{c}+1)s^{b}}\int_{0}^{1}|w|^{p_{c}+1}\varphi (y,s) {\mathrm{d}}y \nonumber \\
&&+C e^{-\frac{p_{c}-1}{32b} s}\hspace{-0,3cm}\int_{0}^{1}\frac{(\partial_{s}w)^2}{1-y^2}y^{N-1}(1-y^2)^{\varepsilon_{1}}{\mathrm{d}}y+C e^{-\frac{p_{c}-1}{16b} s}\hspace{-0,3cm}\int_{0}^{1}(\partial_{y}w)^2y^{N-1}{\mathrm{d}}y+\frac{C}{s^{b}} . \nonumber 
\end{eqnarray}
 Recalling that
 $$L(w(s),s)=\exp\Big(\frac{p_{c}+3}{2(b-1)s^{b-1}}\Big) K(w(s),s)+\frac{\sigma}{s^{b-1}}\,\,\,\,{\rm with }\,\,\,b\in (1,a),$$
a derivative in time of the expression of $L(w(s),s)$ give birth to the following equality:
\begin{eqnarray*}
\frac{d}{ds}(L(w(s),s))&=&-\frac{p_{c}+3}{2s^{b}}\exp\Big(\frac{p_{c}+3}{2(b-1)s^{b-1}}\Big) K(w(s),s)\\
&&+\exp\Big(\frac{p_{c}+3}{2(b-1)s^{b-1}}\Big)\frac{d}{ds}(K(w(s),s))-\frac{\sigma (b-1)}{s^{b}}.
\end{eqnarray*}
Since, for all $s\geq \max(S_{4},s_{0})$, we have $1\leq \exp\Big(\frac{p_{c}+3}{2(b-1)s^{b-1}}\Big)\leq \exp\Big(\frac{p_{c}+3}{2(b-1)}\Big)$ and by exploiting equation $\eqref{ne}$, we can see 
 \begin{eqnarray*}
 \frac{d}{ds}(L(w(s),s))&\leq &-\frac{2(p_{c}-1)}{(p_{c}+15)s^{b}}\int_{0}^{1}(\partial_{s}w)^2\frac{y^{2}\varphi (y,s)}{1-y^2}{\mathrm{d}}y-\frac{p_{c}-1}{16(p_{c}+1)s^{b}}\int_{0}^{1}|w|^{p_{c}+1}\varphi (y,s) {\mathrm{d}}y\\
&&-\frac{p_{c}-1}{16s^{b}}\int_{0}^{1}(\partial_{y}w)^2(1-y^2)\varphi (y,s) {\mathrm{d}}y+C e^{-\frac{p_{c}-1}{16b} s}\int_{0}^{1}(\partial_{y}w)^2y^{N-1}{\mathrm{d}}y\\
&&+C e^{-\frac{p_{c}-1}{32b} s}\int_{0}^{1}\frac{(\partial_{s}w)^2}{1-y^2}y^{N-1}(1-y^2)^{\varepsilon_{1}}{\mathrm{d}}y+\frac{C-\sigma (b-1)}{s^{b}}.
\end{eqnarray*}
We integrate now between $s$ and $s+1$:
\begin{eqnarray}\label{nest}
 L(w(s+1),s+1)-L(w(s),s) &\leq &-\frac{p_{c}-1}{16(p_{c}+1)(2s)^{b}}\int_{s}^{s+1}\int_{0}^{1}|w(y,\tau)|^{p_{c}+1}\varphi (y,\tau) {\mathrm{d}}y{\mathrm{d}}\tau \nonumber \\
&&-\frac{2(p_{c}-1)}{(p_{c}+15)(2s)^{b}}\int_{s}^{s+1}\int_{0}^{1}(\partial_{s}w(y,\tau))^2\frac{y^{2}\varphi (y,\tau)}{1-y^2}{\mathrm{d}}y{\mathrm{d}}\tau \nonumber \\
&&-\frac{p_{c}-1}{16(2s)^{b}}\int_{s}^{s+1}\int_{0}^{1}(\partial_{y}w)^2(1-y^2)\varphi (y,\tau) {\mathrm{d}}y{\mathrm{d}}\tau\\
&&+\underbrace{C e^{-\frac{p_{c}-1}{32b} s}\int_{s}^{s+1}\int_{0}^{1}\frac{(\partial_{s}w)^2}{1-y^2}y^{N-1}(1-y^2)^{\varepsilon_{1}}{\mathrm{d}}y{\mathrm{d}}\tau }_{I_{1}(s)}\nonumber \\
&& \underbrace{+C e^{-\frac{p_{c}-1}{16b} s}\int_{s}^{s+1}\int_{0}^{1}(\partial_{y}w)^2y^{N-1}{\mathrm{d}}y{\mathrm{d}}\tau }_{I_{2}(s)}+\frac{C-\sigma (b-1)}{s^{b}}.\nonumber 
\end{eqnarray}
To end the proof of Proposition 3.1 we need just to estimate $I_{1}(s)$ and $I_{2}(s)$.
\\According to inequality $\eqref{w2}$ in Theorem 1, we choose $\varepsilon_{1}=\frac{p_{c}-1}{(p_{c}+3)32b}\in (0,\frac{1}{2})$, to deduce that:
\begin{equation}\label{nest1}
I_{1}(s)\leq Ce^{-\frac{p_{c}-1}{64b} s} \leq \frac{C}{s^{b}},
\end{equation}
and we choose $\eta =\frac{p_{c}-1}{(p_{c}+3)16b}$ in $\eqref{w3}$ of Theorem 1, to deduce that 
\begin{equation}\label{nest2}
I_{2}(s) \leq Ce^{-\frac{p_{c}-1}{32b} s}\leq \frac{C}{s^{b}}.
\end{equation}
We combine $\eqref{nest}$, $\eqref{nest1}$ and $\eqref{nest2}$ to obtain the following inequality 
\begin{eqnarray*}
 L(w(s+1),s+1)-L(w(s),s) &\leq &-\frac{p_{c}-1}{16(2s)^{b}}\int_{s}^{s+1}\int_{0}^{1}(\partial_{y}w)^2(1-y^2)\varphi (y,\tau) {\mathrm{d}}y{\mathrm{d}}\tau \\
&&-\frac{2(p_{c}-1)}{(p_{c}+15)(2s)^{b}}\int_{s}^{s+1}\hspace{-0,3cm}\int_{0}^{1}(\partial_{s}w(y,\tau))^2\frac{y^{2}\varphi (y,\tau)}{1-y^2}{\mathrm{d}}y{\mathrm{d}}\tau \\
&& -\frac{p_{c}-1}{16(p_{c}+1)(2s)^{b}}\int_{s}^{s+1}\hspace{-0,3cm}\int_{0}^{1}|w(y,\tau)|^{p_{c}+1}\varphi (y,\tau) {\mathrm{d}}y{\mathrm{d}}\tau\\
&&+\frac{C-\sigma (b-1)}{s^{b}},
\end{eqnarray*}
finally, we choose $\lambda_{1}=\min(\frac{p_{c}-1}{16(p_{c}+1)2^{b}},\frac{p_{c}-1}{(p_{c}+15)2^{b-1}},\frac{p_{c}-1}{2^{4+b}})=\frac{p_{c}-1}{16(p_{c}+1)2^{b}}$ and $\sigma $ large enough so that $C-\sigma (b-1) \leq 0$, to deduce that for all $s\geq \max(S_{4},s_{0})$ inequality  $\eqref{new}$ holds. This ends the proof of the first point $\eqref{new}$ of Proposition 3.1.
\\To end the proof of the last point of Proposition 3.1, we refer the reader to \cite{X}. Let us mention that our proof strongly relies on the fact that $p_{c}\equiv 1+\frac{4}{N-1}<1+\frac{4}{N-2}$.
\Box
\subsection{An exponential bound to the time average of the $L^{p_{c}+1}$ norm of $w$ with singular weight}
In this subsection we prove Proposition 3 which allows to prove $(iii)$ of Theorem 4 where we use essentially the Pohozaev identity. To do that, we need to introduce for all $\eta \in (0,1)$ the following new functional $N_{\eta}(w(s))$ defined by:
\begin{equation}\label{L}
N_{\eta}(w(s))=\int_{0}^{1}\Big((y\partial_{y}w)^2+y\partial_{y}w \partial_sw
\Big)\Psi_{\eta}(y) {\mathrm{d}}y,
\end{equation}
with 
\begin{equation}\label{psii}
\Psi_{\eta}(y)=y^{N-1}(1-y^2)^{\eta}.
\end{equation}
We begin by estimating the time derivative of $N_{\eta}(w(s))$ in the following lemma:
\begin{lem}\label{L1}
For all $\eta \in (0,1)$, we have for all $s \geq \max(s_{0},1)$
\begin{eqnarray}\label{p6}
\frac{d}{ds}(N_{\eta}(w(s)))&=& \frac{N-2}{2}\int_{0}^{1}(\partial_{y}  w)^2\Psi_{\eta}(y) {\mathrm{d}}y+(\eta -\frac{N}{2})\int_{0}^{1}(y\partial_{y}  w)^2\Psi_{\eta}(y) {\mathrm{d}}y\\
&& -\frac{N}2\int_{0}^{1}  (\partial_sw)^2 \Psi_{\eta}(y) {\mathrm{d}}y +\eta \int_{0}^{1}(\partial_sw)^2
\frac{y^{2}\Psi_{\eta}(y)}{1-y^2} {\mathrm{d}}y\nonumber\\
&&-\frac{2(p_{c}+1)}{(p_{c}-1)^2}\int_{0}^{1}y\partial_{y}
ww\Psi_{\eta}(y){\mathrm{d}}y -\frac{p_{c}+3}{p_{c}-1} \int_{0}^{1}y\partial_{y} w
\partial_{s}w\Psi_{\eta}(y) {\mathrm{d}}y\nonumber \\
&& -\frac{N}{p_{c}+1}\int_{0}^{1} |w|^{p_{c}+1} \Psi_{\eta}(y) {\mathrm{d}}y
+\frac{2\eta}{p_{c}+1}\int_{0}^{1} |w|^{p_{c}+1}
\frac{y^{2}\Psi_{\eta}(y)}{1-y^2} {\mathrm{d}}y\nonumber\\
 && -Ne^{\frac{-2(p_{c}+1)s}{p_{c}-1}} \int_{0}^{1}F(e^{\frac{2s}{p_{c}-1}}w)\Psi_{\eta}(y)
{\mathrm{d}}y +2\eta e^{\frac{-2(p_{c}+1)s}{p_{c}-1}}
\hspace{-0,2cm}\int_{0}^{1}\hspace{-0,2cm}F(e^{\frac{2s}{p_{c}-1}}w)\frac{y^{2}\Psi_{\eta}(y)}{1-y^2} 
{\mathrm{d}}y,\nonumber
\end{eqnarray}
where $\Psi_{\eta}(y)$ is defined in \eqref{psii}.
\end{lem}

\bigskip

{\it Proof:} Note that  $N_{\eta}(w(s))$ is a differentiable
function for all
 $s\ge s_0$, we have
\begin{eqnarray}\label{p'''}
\frac{d}{ds}(N_{\eta}(w(s)))&=& 2 \int_{0}^{1}\partial_{y}w \partial^2_{y,s} w^2 y^2\Psi_{\eta}(y){\mathrm{d}}y\nonumber\\
&&+\int_{0}^{1}\Big(y\partial_{y} w
\partial^2_sw+y\partial^2_{y,s}w \partial_sw
\Big)\Psi_{\eta}(y) {\mathrm{d}}y.
\end{eqnarray}
Since we see from integration by parts that
$$\int_{0}^{1} y\partial^2_{y,s}w \partial_sw
\Psi_{\eta}(y) {\mathrm{d}}y=-\frac{N}2\int_{0}^{1}  (\partial_sw)^2 \Psi_{\eta}(y) {\mathrm{d}}y
+\eta \int_{0}^{1}  (\partial_sw)^2 \frac{y^{2}\Psi_{\eta}(y)}{1-y^2}
{\mathrm{d}}y.$$
Combining this equality with \eqref{p'''} to write
\begin{eqnarray}\label{3}
\frac{d}{ds}(N_{\eta}(w(s)))&=& -\frac{N}2\int_{0}^{1}  (\partial_sw)^2 \Psi_{\eta}(y) {\mathrm{d}}y+\eta \int_{0}^{1}  (\partial_sw)^2 \frac{y^{2}\Psi_{\eta}(y)}{1-y^2}
{\mathrm{d}}y \nonumber\\
&&+ \int_{0}^{1}y\partial_{y}w(\partial^2_sw +2y\partial^2_{y,s} w)\Psi_{\eta}(y){\mathrm{d}}y.
\end{eqnarray}
By using \eqref{p'} and integrating by parts we have
\begin{eqnarray*}
\frac{d}{ds}(N_{\eta}(w(s)))&=&-\frac{N}2\int_{0}^{1}  (\partial_sw)^2
\Psi_{\eta}(y) {\mathrm{d}}y +\eta \int_{0}^{1}  (\partial_sw)^2
\frac{y^{2}\Psi_{\eta}(y)}{1-y^2} {\mathrm{d}}y\nonumber\\
&&+  \int_{0}^{1}y\partial_{y}  w \partial_{y}(\Psi_{\eta}(y)(1-y^2)\partial_{y}w){\mathrm{d}}y
+2\eta \int_{0}^{1}(y\partial_{y}  w)^2\Psi_{\eta}(y){\mathrm{d}}y\\
&&-\frac{2(p_{c}+1)}{(p_{c}-1)^2}\int_{0}^{1}y\partial_{y}
ww\Psi_{\eta}(y){\mathrm{d}}y
 +\int_{0}^{1}y\partial_{y} w |w|^{p_{c}-1}w
\Psi_{\eta}(y) {\mathrm{d}}y\\
 &&-\frac{p_{c}+3}{p_{c}-1}
\int_{0}^{1}y\partial_{y}  w \partial_{s}w\Psi_{\eta}(y) {\mathrm{d}}y +
e^{\frac{-2p_{c}s}{p_{c}-1}}\int_{0}^{1}y\partial_{y}  w
f(e^{\frac{2s}{p_{c}-1}}w)\Psi_{\eta}(y) {\mathrm{d}}y.
\end{eqnarray*}
Some simple integration by parts ends the proof of Lemma 3.4. 
\Box

\medskip

Now, we are able to deduce Proposition 3.\\

{\it Proof of Proposition 3:} Let $s\ge \max (S_5, s_0)$, $s_3=s_3(s)\in [s-1,s]$  
and \\$s_4=s_4(s)\in [s+1,s+2]$ to be chosen later. From Lemma 3.4 we can see: 
\begin{eqnarray}\label{p4}
\frac{2\eta}{p_{c}+1}\int_{s}^{s+1}\hspace{-0,3cm}\int_{0}^{1} |w|^{p_{c}+1}
\frac{\Psi_{\eta}(y)}{1-y^2} {\mathrm{d}}y {\mathrm{d}}\tau \hspace{-0,3cm} &\le &\hspace{-0,3cm} C\int_{s_3}^{s_4}\!\int_{0}^{1}\Big((y\partial_{y}  w)^2+(\partial_sw)^2+|w|^{p_{c}+1} \Big)\Psi_{\eta}(y) {\mathrm{d}}y{\mathrm{d}}\tau   \nonumber\\
 &&+R_1(s)+C\Big(R_2(s)+R_3(s)+R_4(s)\Big),
\end{eqnarray}
with
\begin{eqnarray*}
R_1(s)&=&N_{\eta}(w(s_4))-N_{\eta}(w(s_3)),\\
R_2(s)&=&\int_{s_3}^{s_4}\int_{0}^{1}y|\partial_{y}
w||w|\Psi_{\eta}(y){\mathrm{d}}y{\mathrm{d}}\tau ,\\
R_3(s)&=& \int_{s_3}^{s_4}\int_{0}^{1}y|\partial_{y} w|
|\partial_{s}w|\Psi_{\eta}(y) {\mathrm{d}}y{\mathrm{d}}\tau ,\\
\displaystyle R_4(s)&=&\int_{s_3}^{s_4}e^{\frac{-2(p_{c}+1)\tau}{p_{c}-1}} \Big(\int_{0}^{1}|F(e^{\frac{2\tau}{p_{c}-1}}w)|\Psi_{\eta}(y)
{\mathrm{d}}y+\int_{0}^{1}|F(e^{\frac{2\tau}{p_{c}-1}}w)|\frac{y^2\Psi_{\eta}(y)}{1-y^2} {\mathrm{d}}y\Big){\mathrm{d}}\tau .
\end{eqnarray*}
Now, we control all the terms on the right-hand side
of the relation (\ref{p4}). Note that, by the expression (\ref{L}) of $N_{\eta}(w(s))$ and using the Cauchy-Schwarz inequality,
we can write
\begin{equation}\label{I1}
-N_{\eta}(w(s_3))\le  \displaystyle\int_{0}^{1}(\partial_s
w(s_3))^2y^{N-1} {\mathrm{d}}y.
\end{equation}
By using the mean value theorem, let us choose  $s_3=s_3(s)\in
[s-1,s]$ such that
\begin{equation}\label{q4}\displaystyle\int_{s-1}^{s}
\displaystyle\int_{0}^{1}(\partial_s w(\tau))^2y^{N-1} 
{\mathrm{d}}y{\mathrm{d}}\tau=\displaystyle\int_{0}^{1}(\partial_s
w(s_3))^2y^{N-1}  {\mathrm{d}}y.
\end{equation}
In view of Theorem 1, (\ref{I1}) and (\ref{q4}) we write
\begin{equation}\label{I1111}
-N_{\eta}(w(s_3))\le Ce^{\eta \frac{p_{c}+3}2s}.
\end{equation}
From (\ref{L}) and the fact that $ab\le a^2+b^2$, we write
\begin{equation}\label{I02}
N_{\eta}(w(s_4))\le
C\displaystyle\int_{0}^{1}\Big((\partial_s w(s_4))^2 +(y\partial_{y} w(s_4))^2\Big) y^{N-1} {\mathrm{d}}y.
\end{equation}
Similarly, by using the mean value theorem, we choose
$s_4=s_4(s)\in [s+1,s+2]$ such that
\begin{equation}\label{q44}
 \noindent \footnotesize{ \int_{0}^{1}\Big((\partial_s w(s_4))^2+(y\partial_{y} w(s_4))^2\Big)
y^{N-1}{\mathrm{d}}y = \int_{s+1}^{s+2}\hspace{-0,2cm}\int_{0}^{1}\Big((\partial_s
w(\tau))^2 
+(y\partial_{y} w(\tau))^2\Big)y^{N-1} {\mathrm{d}}y{\mathrm{d}}\tau.}
\end{equation}
Theorem 1, (\ref{I02}) and (\ref{q44}) implies that
\begin{equation}\label{I2}
N_{\eta}(w(s_4))\le Ce^{\eta \frac{p_{c}+3}2s}.
\end{equation}
By combining (\ref{I1111}) and(\ref{I2}), we deduce that
\begin{equation}\label{p''}
R_1(s)\le Ce^{\eta \frac{p_{c}+3}2s}.
\end{equation}
By the Cauchy-Schwarz and the Young inequality, we can see that 
$$R_2(s)\le C\int_{s_3}^{s_4}\int_{0}^{1}(y\partial_{y}  w)^2\Psi_{\eta}(y) {\mathrm{d}}y{\mathrm{d}}\tau+C\int_{s_3}^{s_4}\int_{0}^{1} |w|^{p_{c}+1} \Psi_{\eta}(y){\mathrm{d}}y{\mathrm{d}}\tau+C.$$
We use again the Cauchy-Schwarz inequality, we obtain:
$$R_3(s)\le C\int_{s_3}^{s_4}\int_{0}^{1}(y\partial_{y}  w)^2\Psi_{\eta}(y) {\mathrm{d}}y{\mathrm{d}}\tau+C\int_{s_3}^{s_4}\int_{0}^{1} (\partial_sw)^2 \Psi_{\eta}(y) {\mathrm{d}}y{\mathrm{d}}\tau.$$
Since $s_3\in [s-1,s]$ and $s_4\in [s+1,s+2]$, from Theorem 1 we obtain
\begin{equation}\label{14}
R_2(s)+R_3(s)\le Ce^{\eta \frac{p_{c}+3}2s}.
\end{equation}
Finally, it remains only to control the  term  $R_4(s)$. Clearly this term verifies the following equality 
$$R_4(s)=\int_{s_3}^{s_4}e^{\frac{-2(p_{c}+1)\tau}{p_{c}-1}}\int_{0}^{1}|F(e^{\frac{2\tau}{p_{c}-1}}w)|\frac{\Psi_{\eta}(y)}{1-y^2} {\mathrm{d}}y{\mathrm{d}}\tau ,$$
Similarly to $\eqref{200}$, we can write
\begin{equation}
R_4(s)\le C\int_{s_3}^{s_4}e^{\frac{-2(p_{c}+1)\tau}{p_{c}-1}}\int_{0}^{1}\frac{\Psi_{\eta}(y)}{1-y^2} {\mathrm{d}}y{\mathrm{d}}\tau+ \frac{C}{s^{a}}\int_{s_3}^{s_4}\int_{0}^{1}|w|^{p_{c}+1}\frac{\Psi_{\eta}(y)}{1-y^2} {\mathrm{d}}y{\mathrm{d}}\tau.
\end{equation}
We can remark that there exists $S_{6}\geq S_{5}$ such that for all $s\geq \max(S_{6},s_{0})$ we have 
\begin{equation}\label{p5}
R_4(s)\le  C+ \frac{\eta}{p_{c}+1}\int_{s_3}^{s_4}\int_{0}^{1}|w|^{p_{c}+1}\frac{\Psi_{\eta}(y)}{1-y^2} {\mathrm{d}}y{\mathrm{d}}\tau.
\end{equation}
Now, we are able to conclude the proof of Proposition 3. By combining (\ref{p4}), (\ref{p''}), (\ref{14}) and (\ref{p5}) with Theorem 1 we get the desired estimate $\eqref{w1}$.
\Box
\subsection{Proof of Theorem 4 }
 We define the following time
\begin{equation} 
t_{1}(0)=max(T(0)-e^{-S_{5}},0).\\
\end{equation}
Since $b\in (1,a)$, according to the Proposition 3.1, we obtain the following corollary which summarizes the principle properties of $K(w(s),s)$.
\begin{cor}{\bf (Estimate on $K(w(s),s)$).}
For all $b\in (1,a)$, there exists 
$t_{1}(0)\in [0,T(0))$ such that, for all $T_{0}\in (t_{1}(0),T(0)]$, for all $s \geq -\log (T_{0}-t_{1}(0))$ and $y\in (0,1)$ we have
$$-C\leq K(w(s),s) \leq \Big(\theta +K(w(s_{0}),s_{0}) \Big)s^{b},$$
$$\int_{s}^{s+1}\int_{0}^{1}(\partial_{s} w(y,\tau))^2\frac{y^{2}\varphi (y,\tau)}{1-y^2}{\mathrm{d}}y{\mathrm{d}}\tau\leq C\Big(\theta +K(w(s_{0}),s_{0} ) \Big)s^{b}, $$
$$ \int_{s}^{s+1}\int_{0}^{\frac{1}{2}}|w(y,\tau)|^{p_{c}+1} {\mathrm{d}}y{\mathrm{d}}\tau +\int_{s}^{s+1}\int_{0}^{\frac{1}{2}}|\partial_{y} w(y,\tau)|^2 {\mathrm{d}}y{\mathrm{d}}\tau\leq C\Big(\theta +K(w(s_{0}),s_{0} ) \Big)s^{b}, $$
where $w=w_{0,T_{0}}$ is defined in $\eqref{0420}$.

\end{cor}
\begin{rem}
Using the definition $\eqref{0420}$ of $w=w_{0,T_{0}}$, we write easily 
$$C\theta +CK(w(s_{0}),s_{0}) \leq K_{0},$$
where $ K_{0}= K_{0}\Big(\eta,T_{0}-t_{1}(0),\|(u(t_{1}(0)),\partial_{t}u(t_{1}(0)))\|_{H^{1}\times L^{2}((-\frac{T_{0}-t_{1}(0)}{\delta_{0}(0)},\frac{T_{0}-t_{1}(0)}{\delta_{0}(0)}))}\Big)$ and $\delta_{0}(0)\in (0,1)$ is defined in $\eqref{17}$.
\end{rem}
{\it Proof of Theorem 4}:
For the deduction of the proof of Theorem 4, we proceed in two steps the first one is devoted to conclude items $(i)$ and $(ii)$ of Theorem 4 which is similar to the deduction of Theorem 1. The second step is devoted to the deduction of item $(iii)$ of Theorem 4 which is different to items $(i)$ and $(ii)$ where we use Proposition 3 and 3.1. 
\paragraph{Proof of items $(i)$ and $(ii)$ of Theorem 4:}
Note that the estimate on the space-time $L^{2}$ norm of $\partial_{s} w$ was already proved in Corollary 3.5. Thus we focus on the space-time $L^{p_{c}+1}$ norm of $w$ and  $L^{2}$ norm of $\partial_{y} w$. This estimate proved in Corollary 3.5 but just for the space-time $L^{p_{c}+1}$ norm of $w$ and  $L^{2}$ norm of $\partial_{y}w$ in $(0,\frac{1}{2})$. To extend this estimate from $(0,\frac{1}{2})$ to $(0,1)$ we refer the reader to Merle and Zaag \cite{kl} (unperturbed case) and Hamza and Zaag \cite{MH1} (perturbed case), where they introduce a new covering argument to extend the estimate of any known space $L^{q}$ norm of $w$, $\partial_{s} w$, or $\partial_{y} w$, from $(0,\frac{1}{2})$ to $(0,1)$.

\paragraph{Proof of item $(iii)$ of Theorem 4:}
Our concern now is to prove $(iii)$ of Theorem 4, to do that, we divide the interval $(0,1)$ into two parts:
$$B_{3}(s)=\{y \in (0,1)\,\,|\,\, 1-y^2\leq  e^{- s}\}\,\,{\rm and }\,\, B_{4}(s)=\{y \in (0,1)\,\,|\,\, 1-y^2>  e^{- s}\}.$$
On the one hand, if $y\in B_{3}(s)$, 
\begin{eqnarray}\label{B4}
\int_{s}^{s+1}\int_{B_{3}(s)}|w|^{p_{c}+1}y^{N-1}{\mathrm{d}}y{\mathrm{d}}\tau &\leq &  e^{- \varepsilon_{2} s}\int_{s}^{s+1}\int_{B_{3}(s)}\frac{|w|^{p_{c}+1}}{( 1-y^2)^{\varepsilon_{2}}} y^{N-1}{\mathrm{d}}y{\mathrm{d}}\tau \nonumber\\
&&\leq e^{- \varepsilon_{2} s}\int_{s}^{s+1}\int_{0}^{1}\frac{|w|^{p_{c}+1}}{( 1-y^2)^{\varepsilon_{2}}} y^{N-1}{\mathrm{d}}y{\mathrm{d}}\tau ,
\end{eqnarray}
for all $\varepsilon_{2}\in (0,1)$.
\\We are now in position to apply Proposition \ref{p}, to get
\begin{equation}\label{B6}
\int_{s}^{s+1}\int_{B_{3}(s)}|w|^{p_{c}+1}y^{N-1}{\mathrm{d}}y{\mathrm{d}}\tau \leq C e^{- \varepsilon_{2} s}e^{(1- \varepsilon_{2} )\frac{p_{c}+3}{2} s},
\end{equation}
now we combine (\ref{B4}), (\ref{B6}) and we choose $\varepsilon_{2} =\frac{p_{c}+3}{p_{c}+5}\in (0,1)$, we deduce that 
\begin{equation}\label{B20}
\int_{s}^{s+1}\int_{B_{3}(s)}|w|^{p_{c}+1}y^{N-1}{\mathrm{d}}y \leq  C.
\end{equation}
On the other hand, if $y\in B_{4}(s)$, by using the following equality 
$$\int_{s}^{s+1}\int_{B_{4}(s)}|w|^{p_{c}+1}y^{N-1}  {\mathrm{d}}y{\mathrm{d}}\tau =\int_{s}^{s+1}\int_{B_{4}(s)}|w|^{p_{c}+1}\frac{\varphi (y,s)}{(1-y^2)^{\frac{1}{s^{b}}}}  {\mathrm{d}}y {\mathrm{d}}\tau ,$$
the fact that, for all $y\in B_{4}(s)$ we have $\frac{1}{(1-y^2)^{\frac{1}{s^{b}}}}\leq C$, and Proposition 3.1 we write
\begin{equation}\label{B21}
\int_{s}^{s+1}\int_{B_{4}(s)}|w|^{p_{c}+1}y^{N-1}  {\mathrm{d}}y{\mathrm{d}}\tau \leq C\int_{s}^{s+1}\int_{0}^{1}|w|^{p_{c}+1}\varphi (y,s) {\mathrm{d}}y {\mathrm{d}}\tau\leq Cs^{b}.
\end{equation}
Inequality (\ref{B20}) and (\ref{B21}) gives rise to $(iii)$ of Theorem 4.
\Box
\section{Proof of Theorem 6}
 This section is devoted to conclude the proof of Theorem 6 when $a> 2$ and $U$ is a radial blow-up solution of $\eqref{oy}$.
\\Firstly, according to the change of variable $\eqref{0420}$ we write equation $\eqref{p'}$ in the following form:
\begin{eqnarray}\label{nn}
\partial^2_{s}w &=&\frac{1}{y^{N-1}} \partial_{y}(y^{N-1}(1-y^2)\partial_{y}w)-\frac{2(p_{c}+1)}{(p_{c}-1)^2}w+|w|^{p_{c}-1}w\nonumber \\
&&-\frac{p_{c}+3}{p_{c}-1}\partial_{s}w-2y \partial^2_{y,s}w+e^{\frac{-2p_{c}s}{p_{c}-1}}f(e^{\frac{2s}{p_{c}-1}}w).
\end{eqnarray}
Secondly, we introduce the following functional:
\begin{equation}\label{kdv1}
H_{0}(w(s),s)=E_{0}(w(s),s)+\frac{1 }{s^{\frac{a-b-1}{2}}},\,\,\,{\rm with }\,\,a>2\,\,\,\,\,\,b\in (1,a)
\end{equation}
and
\begin{eqnarray}\label{nn10}
E_{0}(w(s),s)&=&\hspace{-0,4cm}\int_{0}^{1}\Big(\frac{1}{2}(\partial_{s}w)^2+\frac{1}{2}(\partial_{y}w)^2(1-y^2)+\frac{p_{c}+1}{(p_{c}-1)^2}w^2-\frac{1}{p_{c}+1}|w|^{p_{c}+1}\Big)y^{N-1} {\mathrm{d}}y \nonumber \\
&&-e^{\frac{-2(p_{c}+1)s}{p_{c}-1}} \int_{0}^{1} F(e^{\frac{2s}{p_{c}-1}}w)y^{N-1} {\mathrm{d}}y.
\end{eqnarray}
This section is divided into two parts:
\begin{itemize}
\item Based upon Theorem 4, we prove that $H_{0}(w(s),s)$ is a Lyapunov functional for equation $\eqref{nn}$, which allows to give a blow up criterion for this equation.
\item Finally, we conclude Theorem 6 when $a> 2$ which is the main goal of this paper.
\end{itemize}
\subsection{A Lyapunov functional for equation $\eqref{nn}$}
 We begin this subsection by the following lemma:
\begin{lem}
For all $b\in (1,a)$, we have for all $s \geq -\log (T(0)-t_{1}(0))$
\begin{equation}\label{nn3}
\frac{d}{ds}(E_{0}(w(s),s))= -(\partial_{s}w(1,s))^2+\Sigma_{4}(s),
\end{equation}
with 
$$\Sigma_{4}(s)\leq \frac{C}{s^a}\int_{0}^{1}|w|^{p_{c}+1}y^{N-1}  {\mathrm{d}}y+Ce^{-\frac{p_{c}+1}{p_{c}-1}s}.$$
\end{lem} 
{\it Proof}: Multiplying $\eqref{nn}$ by $\partial_{s}wy^{N-1}$ and we integrate over $(0,1)$, we obtain $\eqref{nn3}$ with
\begin{equation}\label{nn1}
\noindent \footnotesize{ \Sigma_{4}(s)=\frac{2(p_{c}+1)}{p_{c}-1}e^{-\frac{2(p_{c}+1)s}{p_{c}-1}}\int_{0}^{1} \hspace{-0,2cm} F(e^{\frac{2s}{p_{c}-1}}w) y^{N-1}{\mathrm{d}}y-\frac{2e^{-\frac{2p_{c}s}{p_{c}-1}}}{p_{c}-1}\int_{0}^{1}\hspace{-0,2cm}wf(e^{\frac{2s}{p_{c}-1}}w)y^{N-1}{\mathrm{d}}y.}
\end{equation}
Inequality $\eqref{200}$ induces to
\begin{equation}\label{nn2}
 \Sigma_{4}(s)\leq  \frac{C}{s^a}\int_{0}^{1}|w|^{p_{c}+1}y^{N-1}  {\mathrm{d}}y+Ce^{-\frac{p_{c}+1}{p_{c}-1}s}.
\end{equation}
According to $\eqref{nn1}$ and $\eqref{nn2}$ we get the result.
\Box

\medskip

With Lemma 4.1 and Theorem 4 we are in position to prove that $H_{0}(w(s),s)$ is a Lyapunov functional for equation $\eqref{nn}$.
\begin{prop}
For all $b\in (1,a)$, there exists $S_{7}\geq S_{6}$, such that $H_{0}(w(s),s)$ defined in $\eqref{kdv1}$ satisfies for all $s\geq \max(S_{7},-\log (T_{0}-t_{1}(0)))$
$$  H_{0}(w(s+1),s+1)- H_{0}(w(s),s) \leq -\int_{s}^{s+1}(\partial_{s}w(1,\tau))^2{\mathrm{d}}\tau ,$$
where $w= w_{0,T_{0}}$ defined in $\eqref{0420}$.
Moreover, there exists $S_{8}\geq S_{7}$ such that, for all $s\geq max(S_{8},s_{0})$ we have $H_{0}(w(s),s)\geq 0.$
\end{prop}
{\it Proof}:
By using the expression $\eqref{kdv1}$ of $H_{0}(w(s),s)$, we obtain:
\begin{eqnarray}\label{fin1}
 H_{0}(w(s+1),s+1)-H_{0}(w(s),s) &=&  E_{0}(w(s+1),s+1)-E_{0}(w(s),s)\nonumber \\
&& +\frac{1 }{(s+1)^{\frac{a-b-1}{2}}}-\frac{1 }{s^{\frac{a-b-1}{2}}}. 
 \end{eqnarray}
The idea is to apply the Mean Value theorem to the function $x\longmapsto \frac{1}{x^{\frac{a-b-1}{2}}}$, which is a function of class $\C^{\infty}$ ($x\geq 1$), between $s$ and $s+1$, so we can say that there exists a constant $\gamma\in ] 0,1[ $ such that:
$$ \frac{1 }{(s+1)^{\frac{a-b-1}{2}}}-\frac{1 }{s^{\frac{a-b-1}{2}}}=-\frac{a-b-1}{2(s+\gamma )^{\frac{a-b+1}{2}}}.$$
The simple fact that $\gamma\in ] 0,1[ $ and $s\geq 1$, implies:
\begin{equation}\label{fin10}
-\frac{a-b-1}{2(s+\gamma )^{\frac{a-b+1}{2}}}<-\frac{C}{(s+1)^{\frac{a-b+1}{2}}}\leq -\frac{C }{s^{\frac{a-b+1}{2}}}. 
 \end{equation}
The identity $\eqref{fin1}$ and inequality $\eqref{fin10}$ induces to the following inequality:
\begin{equation}\label{fin2}
\displaystyle H_{0}(w(s+1),s+1)-H_{0}(w(s),s) \leq  E_{0}(w(s+1),s+1)-E_{0}(w(s),s)-\frac{C }{s^{\frac{a-b+1}{2}}}.
\end{equation}
From Lemma 4.1 and inequality $\eqref{fin2}$, we get 
\begin{eqnarray*}
 H_{0}(w(s+1),s+1)-H_{0}(w(s),s) &\leq & -\int_{s}^{s+1}(\partial_{s}w(1,\tau))^2{\mathrm{d}}\tau  +Ce^{-\frac{p_{c}+1}{p_{c}-1}s}\\
 &&+\frac{C}{s^{a}}\int_{s}^{s+1}\int_{0}^{1}|w|^{p_{c}+1}y^{N-1}  {\mathrm{d}}y{\mathrm{d}}\tau -\frac{C }{s^{\frac{a-b+1}{2}}}. 
 \end{eqnarray*}
Besides, from item $(iii)$ of Theorem 4 and the fact that when $s$ is large we have 
\\$e^{-\frac{p_{c}+1}{p_{c}-1}s}\leq \frac{C}{s^{a-b}}$, we write 
$$H_{0}(w(s+1),s+1)-H_{0}(w(s),s) \leq  -\int_{s}^{s+1}(\partial_{s}w(1,\tau))^2{\mathrm{d}}\tau +\frac{C}{s^{a-b}}-\frac{C }{s^{\frac{a-b+1}{2}}}.$$
As we mentionned above, the fact that we fix $a>2$ we can choose $b=\frac{a}{2}>1$ then write:
\begin{eqnarray*}
H_{0}(w(s+1),s+1)-H_{0}(w(s),s) &\leq & -\int_{s}^{s+1}(\partial_{s}w(1,\tau))^2{\mathrm{d}}\tau +\frac{C}{s^{\frac{a}{2}}}-\frac{C }{s^{\frac{a+2}{4}}}\\
 &&=-\int_{s}^{s+1}(\partial_{s}w(1,\tau))^2{\mathrm{d}}\tau +\Big(\frac{C}{s^{\frac{a-2}{4}}}-C\Big)\frac{1}{s^{\frac{a+2}{4}}}.
 \end{eqnarray*}
The fact that $\frac{a-2}{4}>0$ allows us to choose $S_{7}\geq S_{6}$ such that we have $\frac{C}{s^{\frac{a-2}{4}}}-C\leq 0$ to get Proposition 4.2.
To end the proof of the last point of Proposition 4.2, we refer the reader to \cite{cf}. 
\Box
 
 \subsection{Boundedness of the solution in similarity variables}
We prove Theorem 6 here for $a>2$.\,\,Note that the lower bound follows from the finite speed of propagation and the wellposedness in $H^1\times L^2$. For a detailed argument in the similar case of equation $\eqref{1.6}$, (see Lemma 3.1 p.1136 in \cite{kl}).

 \bigskip
 We define the following time 
 $$ t_{2}(0)=\max(T(0)-e^{-S_{8}},0).$$
For some $T_{0}\in (t_{2}(0),T(0)]$, for all $r\in \R^+ $ is such that $r\leq \frac{T_{0}-t_{2}(0)}{\delta_{0}(0)}$, where $\delta_{0}$ is defined in $\eqref{17}$, then we write $w$ instead of $w_{0,T_{0}}$ defined in $\eqref{nn}$. We aim at bounding $\|(w,\partial_{s}w)(s)\|_{H^{1}\times L^{2}}$ for $s$ large.

\begin{cor}
 For all $s \geq -\log(T_{0}-t_{2}(0))$,
 it holds that
\begin{eqnarray*}
-C \leq E_{0}(w(s),s)&\leq &K,\\
 \int_{s}^{s+1}(\partial_{s}w(1,\tau))^2{\mathrm{d}}\tau &\leq &K,\\
\int_{s}^{s+1}\int_{0}^{1} \Big(\partial_{s}w(y,\tau )-\lambda (\tau ,s)w(y,\tau )\Big)^2y^{N-1}{\mathrm{d}}y{\mathrm{d}}\tau &\leq &K,
\end{eqnarray*}
where $0\leq\lambda (\tau ,s)\leq C$, $K=K(T_{0},\|(u(t_{2}(0)),\partial_{t} u(t_{2}(0))\|_{H^1\times L^2((-\frac{T_{0}-t_{2}(0)}{\delta_{0}(0)},\frac{T_{0}-t_{2}(0)}{\delta_{0}(0)}))})$, $C>0$ and $\delta_{0}(0)$ is defined in $\eqref{17}$.
\end{cor}
{\it Proof}: As in \cite{MH1}, from Proposition 4.2 we get the first and the second inequality.
\\For the proof of the last inequality the argument is the same as in the corresponding part, (see Proposition 4.2 p.1147 in \cite{kl}).
\Box

\medskip
 
The proof of Theorem 6 is similar to the one in the unperturbed case treated by Merle and Zaag in \cite{fh3} and \cite{fh4} and also used by Hamza and Zaag in \cite{MH}, \cite{MH1} and Hamza and Saidi \cite{X}.\,\,To be accurate and concise in our results, there is an analogy between the exponential smallness exploited in \cite{MH1} by Hamza and Zaag and the polynomial smallness used here, the unique difference lays in the treatement of the perturbed term which is treated by Hamza and Saidi \cite{X}. Which close the proof of Theorem 6.
 \Box

\appendix
\section{The Hardy-Sobolev inequality}
In this part we are going to prove the following inequality
\begin{equation}\label{A}
 \int_{B}h^2\frac{|y|^2\rho_{\eta}}{1-|y|^2}{\mathrm{d}}y\leq \frac{1}{\eta^2}\int_{B}|\nabla h|^2(1-|y|^2)\rho_{\eta} {\mathrm{d}}y+\frac{N}{\eta}\int_{B} h^2\rho_{\eta} {\mathrm{d}}y.
\end{equation}
{\it Proof}:
Some computations give 
$$y.\nabla \rho_{\eta}=-2\eta \frac{|y|^{2}\rho_{\eta}}{1-|y|^2},$$
which implies that
$$ \int_{B}h^2\frac{|y|^{2}\rho_{\eta}}{1-|y|^2}{\mathrm{d}}y=-\frac{1}{2\eta}\int_{B}h^2y.\nabla \rho_{\eta}{\mathrm{d}}y.$$
If we integrate by part we see:
$$\int_{B}h^2\frac{|y|^{2}\rho_{\eta}}{1-|y|^2}{\mathrm{d}}y=\frac{1}{2\eta}\Big(N\int_{B}h^2\rho_{\eta}{\mathrm{d}}y+2\int_{B}h\nabla h .y\rho_{\eta}{\mathrm{d}}y\Big).$$
According to the Cauchy-Schwarz inequality 
\begin{eqnarray*}
\Big|\int_{B}h\nabla h .y\rho_{\eta} {\mathrm{d}}y\Big|&\leq &\Big(\int_{B}|\nabla h |^2\rho_{\eta}(1-|y|^2){\mathrm{d}}y\Big)^{\frac{1}{2}} \Big(\int_{B}h^2\frac{|y|^{2}\rho_{\eta}}{1-|y|^2}{\mathrm{d}}y\Big)^{\frac{1}{2}}\\
&\leq & \frac{1}{\varepsilon}\int_{B}|\nabla h |^2\rho_{\eta}(1-|y|^2){\mathrm{d}}y+\varepsilon \int_{B}h^2\frac{|y|^{2}\rho_{\eta}}{1-|y|^2}{\mathrm{d}}y,
\end{eqnarray*}
for any $\varepsilon >0$. We can deduce that
$$ \int_{B}h^2\frac{|y|^{2}\rho_{\eta}}{1-|y|^2}{\mathrm{d}}y\leq \frac{1}{2\eta}\Big(\frac{2}{\varepsilon}\int_{B}|\nabla h |^2\rho_{\eta}(1-|y|^2){\mathrm{d}}y+2\varepsilon \int_{B}h^2\frac{|y|^{2}\rho_{\eta}}{1-|y|^2}{\mathrm{d}}y+N\int_{B}h^2\rho_{\eta}{\mathrm{d}}y\Big).$$
Finally if we choose $\varepsilon =\frac{\eta}{2},$ we show that $\eqref{A}$ holds.
\Box

\noindent{\bf Address}:\\
Universit\'e de Tunis El Manar, Facult\'e des Sciences de Tunis, LR03ES04 \'Equations aux d\'eriv\'ees partielles et applications, 2092 Tunis, Tunisie\\
\vspace{-7mm}
\begin{verbatim}
e-mail: ma.hamza@fst.rnu.tn
e-mail: saidi.omar@hotmail.fr
\end{verbatim}

 %%%%%%%%%%%%%%%%%%%%%%%%%%%%%%%%%%%%%%%%%%%%%%%%%%%%%%%%%%%%%%%%%%%%%%%%%%%%%%%%%%%%%%%%%%%%%%%%%%%%%%%%%%%%%%%

 \end{document}